\documentclass[11pt]{article}

\usepackage{amsmath}
\usepackage{amssymb}
\usepackage{amsthm}
\usepackage{mathtools}
\usepackage{tabu}
\usepackage{textcomp}
\usepackage{listings}
\usepackage{graphicx}
\usepackage{tikz}
\usepackage{arydshln}
\usepackage[margin=1in]{geometry}
\usepackage{indentfirst}
\usepackage{array}
\usepackage{xltabular}
\usepackage{caption}
\usepackage{empheq}
\usepackage{footmisc}

\usepackage{chessfss}
\setboardfontsize{15pt}

\definecolor{colori}{RGB}{0,60,0}
\definecolor{colorii}{RGB}{0,120,0}
\definecolor{coloriii}{RGB}{0,180,0}
\definecolor{coloriv}{RGB}{0,240,0}

\newcommand{\varT}{\mathcal{T}}

\newcommand{\tikzpic}[2]{
\begin{tikzpicture}[scale=#1]
    #2
\end{tikzpicture}
}

\newcommand{\wvert}[2]{\filldraw[fill=white](#1,#2)circle(0.1)}
\newcommand{\bvert}[2]{\filldraw[fill=black](#1,#2)circle(0.1)}
\newcommand{\bking}[2]{\node at(#1,#2){\BlackKingOnWhite}}

\newcommand{\drawwtile}[6]{
    \draw (#1,#2) rectangle (#1+1,#2+1);
    \node at (#1+0.2,#2+0.5) [left] {#3};
    \node at (#1+0.5,#2+0.8) [above] {\phantom{$d$} #4 \phantom{$p$}};
    \node at (#1+0.5,#2+0.2) [below] {\phantom{$d$} #5 \phantom{$p$}};
    \node at (#1+0.8,#2+0.5) [right] {#6}
}

\newcommand{\wtileUo}[4]{\hspace{-1em}
    \raisebox{-0.3cm-1em}{
        \tikzpic{0.6}{
            \drawwtile{0}{0}{$0$}{$00$}{$00$}{$0$};
            \wvert{0}{0.5};\wvert{1}{0.5};
        }
    }\hspace{-1em}
}

\newcommand{\wtileUv}[4]{\hspace{-1em}
    \raisebox{-0.3cm-1em}{
        \tikzpic{0.6}{
            \drawwtile{0}{0}{$0$}{$01$}{$00$}{$1$};
            \draw[line width=1pt](1,1)--(1,0.5);
            \wvert{0}{0.5};\bvert{1}{0.5};
        }
    }\hspace{-1em}
}

\newcommand{\wtileo}[4]{\hspace{-1em}
    \raisebox{-0.3cm-1em}{
        \tikzpic{0.6}{
            \drawwtile{0}{0}{$#1$}{$#2$}{$#3$}{$#4$};
            \wvert{0.5}{0.5};
        }
    }\hspace{-1em}
}

\newcommand{\wtilei}[4]{\hspace{-1em}
    \raisebox{-0.3cm-1em}{
        \tikzpic{0.6}{
            \drawwtile{0}{0}{$#1$}{$#2$}{$#3$}{$#4$};
            \draw[line width=1pt](0.5,0.5)--(0,0.5);
            \wvert{0.5}{0.5};
        }
    }\hspace{-1em}
}

\newcommand{\wtileii}[4]{\hspace{-1em}
    \raisebox{-0.3cm-1em}{
        \tikzpic{0.6}{
            \drawwtile{0}{0}{$#1$}{$#2$}{$#3$}{$#4$};
            \draw[line width=1pt](1,0.5)--(0.5,0.5)--(0.5,1);
            \bvert{0.5}{0.5};
        }
    }\hspace{-1em}
}

\newcommand{\wtileiii}[4]{\hspace{-1em}
    \raisebox{-0.3cm-1em}{
        \tikzpic{0.6}{
            \drawwtile{0}{0}{$#1$}{$#2$}{$#3$}{$#4$};
            \draw[line width=1pt](0.5,0)--(0.5,0.5)--(0,0.5);
            \wvert{0.5}{0.5};
        }
    }\hspace{-1em}
}

\newcommand{\wtileiv}[4]{\hspace{-1em}
    \raisebox{-0.3cm-1em}{
        \tikzpic{0.6}{
            \drawwtile{0}{0}{$#1$}{$#2$}{$#3$}{$#4$};
            \draw[line width=1pt](0.5,0)--(0.5,0.5);
            \wvert{0.5}{0.5};
        }
    }\hspace{-1em}
}

\newcommand{\wtilev}[4]{\hspace{-1em}
    \raisebox{-0.3cm-1em}{
        \tikzpic{0.6}{
            \drawwtile{0}{0}{$#1$}{$#2$}{$#3$}{$#4$};
            \draw[line width=1pt](1,0.5)--(0.5,0.5)--(0.5,0);
            \wvert{0.5}{0.5};
        }
    }\hspace{-1em}
}

\newcommand{\gridi}{
    \begin{tikzpicture}[xscale=-0.6,yscale=0.6]
        \draw(0,0)grid(5,3);
        \wvert{0}{3};\bvert{1}{3};\wvert{2}{3};\bvert{3}{3};\wvert{4}{3};\bvert{5}{3};
        \wvert{0}{2};\wvert{1}{2};\bvert{2}{2};\wvert{3}{2};\bvert{4}{2};\wvert{5}{2};
        \bvert{0}{1};\wvert{1}{1};\wvert{2}{1};\bvert{3}{1};\wvert{4}{1};\wvert{5}{1};
        \wvert{0}{0};\bvert{1}{0};\wvert{2}{0};\wvert{3}{0};\bvert{4}{0};\wvert{5}{0};
    \end{tikzpicture}
}

\newcommand{\squgridi}{
    \begin{tikzpicture}[xscale=-0.6,yscale=0.6]
        \draw[help lines,shift={(-0.5,-0.5)}](0,0)grid(6,4);\draw(0,0)grid(5,3);
        \wvert{0}{3};\bvert{1}{3};\wvert{2}{3};\bvert{3}{3};\wvert{4}{3};\bvert{5}{3};
        \wvert{0}{2};\wvert{1}{2};\bvert{2}{2};\wvert{3}{2};\bvert{4}{2};\wvert{5}{2};
        \bvert{0}{1};\wvert{1}{1};\wvert{2}{1};\bvert{3}{1};\wvert{4}{1};\wvert{5}{1};
        \wvert{0}{0};\bvert{1}{0};\wvert{2}{0};\wvert{3}{0};\bvert{4}{0};\wvert{5}{0};
    \end{tikzpicture}
}

\newcommand{\cylgridi}{
    \begin{tikzpicture}[xscale=0.6,yscale=0.6]
        \draw(0,0)grid(5,3);
        \draw(0,3)arc(270:90:0.5 and 0.25);\draw(5,3)arc(-90:90:0.5 and 0.25);\draw(0,3.5)--(5,3.5);
        \draw[gray](0,2)arc(270:90:0.5 and 0.25);\draw[gray](5,2)arc(-90:90:0.5 and 0.25);\draw[gray](0,2.5)--(5,2.5);
        \draw[gray](0,1)arc(270:90:0.5 and 0.25);\draw[gray](5,1)arc(-90:90:0.5 and 0.25);\draw[gray](0,1.5)--(5,1.5);
        \draw[gray](0,0)arc(270:90:0.5 and 0.25);\draw[gray](5,0)arc(-90:90:0.5 and 0.25);\draw[gray](0,0.5)--(5,0.5);

        \wvert{0}{3};\bvert{1}{3};\wvert{2}{3};\bvert{3}{3};\wvert{4}{3};\bvert{5}{3};
        \bvert{0}{2};\wvert{1}{2};\bvert{2}{2};\wvert{3}{2};\bvert{4}{2};\wvert{5}{2};
        \wvert{0}{1};\wvert{1}{1};\wvert{2}{1};\bvert{3}{1};\wvert{4}{1};\bvert{5}{1};
        \wvert{0}{0};\bvert{1}{0};\wvert{2}{0};\wvert{3}{0};\bvert{4}{0};\wvert{5}{0};
    \end{tikzpicture}
}

\newcommand{\gridii}{
    \begin{tikzpicture}[xscale=-0.6,yscale=0.6]
        \draw(0,0)grid(5,3);
        \draw[line width=1pt](-1,1)--(0,1)--(0,2);
        \draw[line width=1pt](0,3)--(1,3)--(1,4);
        \draw[line width=1pt](0,0)--(1,0)--(1,1);
        \draw[line width=1pt](1,2)--(2,2)--(2,3)--(3,3)--(3,4);
        \draw[line width=1pt](2,1)--(3,1)--(3,2)--(4,2)--(4,3)--(5,3)--(5,4);
        \draw[line width=1pt](3,0)--(4,0)--(4,1);

        \wvert{0}{3};\bvert{1}{3};\wvert{2}{3};\bvert{3}{3};\wvert{4}{3};\bvert{5}{3};
        \wvert{0}{2};\wvert{1}{2};\bvert{2}{2};\wvert{3}{2};\bvert{4}{2};\wvert{5}{2};
        \bvert{0}{1};\wvert{1}{1};\wvert{2}{1};\bvert{3}{1};\wvert{4}{1};\wvert{5}{1};
        \wvert{0}{0};\bvert{1}{0};\wvert{2}{0};\wvert{3}{0};\bvert{4}{0};\wvert{5}{0};
    \end{tikzpicture}
}

\newcommand{\kingi}{
    \begin{tikzpicture}[xscale=0.6,yscale=0.6]
        \draw(0,0)grid(5,3);
        \draw(1,0)--(0,1);\draw(2,0)--(0,2);\draw(3,0)--(0,3);\draw(4,0)--(1,3);\draw(5,0)--(2,3);\draw(5,1)--(3,3);\draw(5,2)--(4,3);
        \draw(1,3)--(0,2);\draw(2,3)--(0,1);\draw(3,3)--(0,0);\draw(4,3)--(1,0);\draw(5,3)--(2,0);\draw(5,2)--(3,0);\draw(5,1)--(4,0);
        \wvert{0}{3};\wvert{1}{3};\wvert{2}{3};\bking{3}{3};\wvert{4}{3};\wvert{5}{3};
        \bking{0}{2};\wvert{1}{2};\wvert{2}{2};\wvert{3}{2};\wvert{4}{2};\bking{5}{2};
        \wvert{0}{1};\wvert{1}{1};\wvert{2}{1};\wvert{3}{1};\wvert{4}{1};\wvert{5}{1};
        \wvert{0}{0};\bking{1}{0};\wvert{2}{0};\wvert{3}{0};\bking{4}{0};\wvert{5}{0};
    \end{tikzpicture}
}

\newcommand{\kingii}{
    \begin{tikzpicture}[xscale=0.6,yscale=0.6]
        \draw(0,0)grid(5,3);
        \draw(1,0)--(0,1);\draw(2,0)--(0,2);\draw(3,0)--(0,3);\draw(4,0)--(1,3);\draw(5,0)--(2,3);\draw(5,1)--(3,3);\draw(5,2)--(4,3);
        \draw(1,3)--(0,2);\draw(2,3)--(0,1);\draw(3,3)--(0,0);\draw(4,3)--(1,0);\draw(5,3)--(2,0);\draw(5,2)--(3,0);\draw(5,1)--(4,0);
        \draw[line width=1pt](1,2)--(0,2)--(0,3)--(1,3);
        \draw[line width=1pt](2,0)--(1,0)--(1,1)--(2,1);
        \draw[line width=1pt](4,3)--(3,3)--(3,4)--(4,4);
        \draw[line width=1pt](5,0)--(4,0)--(4,1)--(5,1);
        \draw[line width=1pt](6,2)--(5,2)--(5,3)--(6,3);

        \wvert{0}{3};\wvert{1}{3};\wvert{2}{3};\bking{3}{3};\wvert{4}{3};\wvert{5}{3};
        \bking{0}{2};\wvert{1}{2};\wvert{2}{2};\wvert{3}{2};\wvert{4}{2};\bking{5}{2};
        \wvert{0}{1};\wvert{1}{1};\wvert{2}{1};\wvert{3}{1};\wvert{4}{1};\wvert{5}{1};
        \wvert{0}{0};\bking{1}{0};\wvert{2}{0};\wvert{3}{0};\bking{4}{0};\wvert{5}{0};
    \end{tikzpicture}
}

\newcommand{\hex}[2]{
    \draw[help lines](#1-0.5,#2-0.2667)--(#1,#2-0.5333)--(#1+0.5,#2-0.2667)--(#1+0.5,#2+0.2667)--(#1,#2+0.5333)--(#1-0.5,#2+0.2667)--(#1-0.5,#2-0.2667)
}

\newcommand{\hexbvert}[2]{\hex{#1}{#2};\bvert{#1}{#2}}
\newcommand{\hexwvert}[2]{\hex{#1}{#2};\wvert{#1}{#2}}

\newcommand{\trigridi}{
    \begin{tikzpicture}[xscale=-0.6,yscale=0.6]
        \draw(0.0,0.0)--(5.0,0.0);\draw(0.5,0.8)--(4.5,0.8);\draw(1.0,1.6)--(4.0,1.6);\draw(1.5,2.4)--(3.5,2.4);\draw(2.0,3.2)--(3.0,3.2);
        \draw(0.0,0.0)--(2.5,4.0);\draw(1.0,0.0)--(3.0,3.2);\draw(2.0,0.0)--(3.5,2.4);\draw(3.0,0.0)--(4.0,1.6);\draw(4.0,0.0)--(4.5,0.8);
        \draw(1.0,0.0)--(0.5,0.8);\draw(2.0,0.0)--(1.0,1.6);\draw(3.0,0.0)--(1.5,2.4);\draw(4.0,0.0)--(2.0,3.2);\draw(5.0,0.0)--(2.5,4.0);
                            \bvert{2.5}{4.0};
                        \wvert{2.0}{3.2};\wvert{3.0}{3.2};
                    \wvert{1.5}{2.4};\bvert{2.5}{2.4};\wvert{3.5}{2.4};
                \bvert{1.0}{1.6};\wvert{2.0}{1.6};\wvert{3.0}{1.6};\wvert{4.0}{1.6};
            \wvert{0.5}{0.8};\wvert{1.5}{0.8};\wvert{2.5}{0.8};\bvert{3.5}{0.8};\wvert{4.5}{0.8};
        \wvert{0.0}{0.0};\bvert{1.0}{0.0};\wvert{2.0}{0.0};\wvert{3.0}{0.0};\wvert{4.0}{0.0};\bvert{5.0}{0.0};
    \end{tikzpicture}
}

\newcommand{\hextrigridi}{
    \begin{tikzpicture}[xscale=-0.6,yscale=0.6]
        \draw(0.0,0.0)--(5.0,0.0);\draw(0.5,0.8)--(4.5,0.8);\draw(1.0,1.6)--(4.0,1.6);\draw(1.5,2.4)--(3.5,2.4);\draw(2.0,3.2)--(3.0,3.2);
        \draw(0.0,0.0)--(2.5,4.0);\draw(1.0,0.0)--(3.0,3.2);\draw(2.0,0.0)--(3.5,2.4);\draw(3.0,0.0)--(4.0,1.6);\draw(4.0,0.0)--(4.5,0.8);
        \draw(1.0,0.0)--(0.5,0.8);\draw(2.0,0.0)--(1.0,1.6);\draw(3.0,0.0)--(1.5,2.4);\draw(4.0,0.0)--(2.0,3.2);\draw(5.0,0.0)--(2.5,4.0);
                            \hexbvert{2.5}{4.0};
                        \hexwvert{2.0}{3.2};\hexwvert{3.0}{3.2};
                    \hexwvert{1.5}{2.4};\hexbvert{2.5}{2.4};\hexwvert{3.5}{2.4};
                \hexbvert{1.0}{1.6};\hexwvert{2.0}{1.6};\hexwvert{3.0}{1.6};\hexwvert{4.0}{1.6};
            \hexwvert{0.5}{0.8};\hexwvert{1.5}{0.8};\hexwvert{2.5}{0.8};\hexbvert{3.5}{0.8};\hexwvert{4.5}{0.8};
        \hexwvert{0.0}{0.0};\hexbvert{1.0}{0.0};\hexwvert{2.0}{0.0};\hexwvert{3.0}{0.0};\hexwvert{4.0}{0.0};\hexbvert{5.0}{0.0};
    \end{tikzpicture}
}

\newcommand{\trigridii}{
    \begin{tikzpicture}[xscale=0.6,yscale=0.6]
        \draw(0.0,0.0)--(5.0,0.0);\draw(0.5,0.8)--(5.5,0.8);\draw(1.0,1.6)--(6.0,1.6);\draw(1.5,2.4)--(6.5,2.4);\draw(2.0,3.2)--(7.0,3.2);\draw(2.5,4.0)--(7.5,4.0);
        \draw(0.0,0.0)--(2.5,4.0);\draw(1.0,0.0)--(3.5,4.0);\draw(2.0,0.0)--(4.5,4.0);\draw(3.0,0.0)--(5.5,4.0);\draw(4.0,0.0)--(6.5,4.0);\draw(5.0,0.0)--(7.5,4.0);
        \draw(1.0,0.0)--(0.5,0.8);\draw(2.0,0.0)--(1.0,1.6);\draw(3.0,0.0)--(1.5,2.4);\draw(4.0,0.0)--(2.0,3.2);
        \draw(5.0,0.0)--(2.5,4.0);\draw(5.5,0.8)--(3.5,4.0);\draw(6.0,1.6)--(4.5,4.0);\draw(6.5,2.4)--(5.5,4.0);\draw(7.0,3.2)--(6.5,4.0);

                            \bvert{2.5}{4.0};\wvert{3.5}{4.0};\wvert{4.5}{4.0};\bvert{5.5}{4.0};\wvert{6.5}{4.0};\wvert{7.5}{4.0};
                        \wvert{2.0}{3.2};\wvert{3.0}{3.2};\wvert{4.0}{3.2};\wvert{5.0}{3.2};\wvert{6.0}{3.2};\wvert{7.0}{3.2};
                    \wvert{1.5}{2.4};\bvert{2.5}{2.4};\wvert{3.5}{2.4};\bvert{4.5}{2.4};\wvert{5.5}{2.4};\bvert{6.5}{2.4};
                \bvert{1.0}{1.6};\wvert{2.0}{1.6};\wvert{3.0}{1.6};\wvert{4.0}{1.6};\wvert{5.0}{1.6};\wvert{6.0}{1.6};
            \wvert{0.5}{0.8};\wvert{1.5}{0.8};\wvert{2.5}{0.8};\bvert{3.5}{0.8};\wvert{4.5}{0.8};\wvert{5.5}{0.8};
        \wvert{0.0}{0.0};\bvert{1.0}{0.0};\wvert{2.0}{0.0};\wvert{3.0}{0.0};\wvert{4.0}{0.0};\bvert{5.0}{0.0};
    \end{tikzpicture}
}

\newcommand{\trigridiii}{
    \begin{tikzpicture}[xscale=-0.6,yscale=0.6]
        \draw[help lines](0,0)grid(5,5);

        \draw(0,0)--(5,0);\draw(1,1)--(5,1);\draw(2,2)--(5,2);\draw(3,3)--(5,3);\draw(4,4)--(5,4);
        \draw(0,0)--(5,5);\draw(1,0)--(5,4);\draw(2,0)--(5,3);\draw(3,0)--(5,2);\draw(4,0)--(5,1);
        \draw(1,0)--(1,1);\draw(2,0)--(2,2);\draw(3,0)--(3,3);\draw(4,0)--(4,4);\draw(5,0)--(5,5);

        \draw[line width=1pt](4,5)--(5,5)--(5,6);
        \draw[line width=1pt](3,3)--(4,3)--(4,4);
        \draw[line width=1pt](1,2)--(2,2)--(2,3);
        \draw[line width=1pt](3,1)--(4,1)--(4,2);
        \draw[line width=1pt](0,0)--(1,0)--(1,1);
        \draw[line width=1pt](4,0)--(5,0)--(5,1);

                                                                         \bvert{5}{5};
                                                            \wvert{4}{4};\wvert{5}{4};
                                               \wvert{3}{3};\bvert{4}{3};\wvert{5}{3};
                                  \bvert{2}{2};\wvert{3}{2};\wvert{4}{2};\wvert{5}{2};
                     \wvert{1}{1};\wvert{2}{1};\wvert{3}{1};\bvert{4}{1};\wvert{5}{1};
        \wvert{0}{0};\bvert{1}{0};\wvert{2}{0};\wvert{3}{0};\wvert{4}{0};\bvert{5}{0};
    \end{tikzpicture}
}

\newcommand{\tensorT}[2]{
    \filldraw[fill=gray!20](#1-0.25,#2-0.25)rectangle(#1+0.25,#2+0.25);
    \node at(#1,#2) {$\varT$}
}

\newcommand{\tensorTD}[2]{
    \filldraw[fill=gray!20](#1-0.25,#2-0.25)rectangle(#1+0.25,#2+0.25);
    \node at(#1,#2) {$\varT_\varDelta$}
}

\newcommand{\tensorvoi}[2]{
    \filldraw[fill=gray!20](#1,#2)circle(0.25);
    \filldraw[fill=gray!20, draw=gray!20](#1-0.25,#2-0.25)rectangle(#1+0.25,#2);
    \draw(#1-0.25,#2)--(#1-0.25,#2-0.25)--(#1+0.25,#2-0.25)--(#1+0.25,#2);
    \node at(#1,#2) {$_{0\!+\!1}$}
}

\newcommand{\tensorhoi}[2]{
    \filldraw[fill=gray!20](#1,#2)circle(0.25);
    \filldraw[fill=gray!20, draw=gray!20](#1,#2-0.25)rectangle(#1+0.25,#2+0.25);
    \draw(#1,#2+0.25)--(#1+0.25,#2+0.25)--(#1+0.25,#2-0.25)--(#1,#2-0.25);
    \node at(#1,#2) {$0$}
}

\newcommand{\tensorho}[2]{
    \filldraw[fill=gray!20](#1,#2)circle(0.25);
    \filldraw[fill=gray!20, draw=gray!20](#1-0.25,#2-0.25)rectangle(#1,#2+0.25);
    \draw(#1,#2+0.25)--(#1-0.25,#2+0.25)--(#1-0.25,#2-0.25)--(#1,#2-0.25);
    \node at(#1,#2) {$_{0\!+\!1}$}
}

\newcommand{\tensorvo}[2]{
    \filldraw[fill=gray!20](#1,#2)circle(0.25);
    \filldraw[fill=gray!20, draw=gray!20](#1-0.25,#2)rectangle(#1+0.25,#2+0.25);
    \draw(#1-0.25,#2)--(#1-0.25,#2+0.25)--(#1+0.25,#2+0.25)--(#1+0.25,#2);
    \node at(#1,#2) {$0$}
}

\newcommand{\TN}{
    \begin{tikzpicture}
        \draw(0.1,2.5)grid(2.5,4.9);\draw(3.5,2.5)grid(4.9,4.9);
        \draw(0.1,0.1)grid(2.5,1.5);\draw(3.5,0.1)grid(4.9,1.5);
        \node at(3,1){$\ldots$};\node at(3,3){$\ldots$};\node at(3,4){$\ldots$};
        \node at(1,2){$\vdots$};\node at(2,2){$\vdots$};\node at(4,2){$\vdots$};
        \tensorvo{1}{0};\tensorvo{2}{0};\tensorvo{4}{0};
        \tensorvoi{1}{5};\tensorvoi{2}{5};\tensorvoi{4}{5};
        \tensorhoi{0}{1};\tensorhoi{0}{3};\tensorhoi{0}{4};
        \tensorho{5}{1};\tensorho{5}{3};\tensorho{5}{4};
        \tensorT{1}{4};\tensorT{2}{4};\tensorT{4}{4};
        \tensorT{1}{3};\tensorT{2}{3};\tensorT{4}{3};
        \tensorT{1}{1};\tensorT{2}{1};\tensorT{4}{1};
        \draw[|<->|](0.5,5.5)--(4.5,5.5);\node at(2.5,5.5)[above]{$m$ columns};
        \draw[|<->|](5.5,0.5)--(5.5,4.5);\node at(5.5,2.5)[right]{$n$ rows};
    \end{tikzpicture}
}

\newcommand{\TNtri}{
    \begin{tikzpicture}
        \draw(0.1,2.5)grid(2.9,4.9);
        \draw(0.1,0.1)grid(2.5,1.5);\draw(3.5,0.1)grid(4.9,1.9);
        \filldraw[draw=white,fill=white](2,4) circle (0.05);\draw(2,3.9)--(2,4.1);
        \node at(1,2){$\vdots$};\node at(2,2){$\vdots$};\node at(1,2){$\vdots$};\node at(3,1){$\ldots$};
        \tensorvo{1}{0};\tensorvo{2}{0};\tensorvo{4}{0};
        \tensorvoi{1}{5};\tensorvoi{2}{5};\tensorvoi{4}{2};
        \tensorhoi{0}{1};\tensorhoi{0}{3};\tensorhoi{0}{4};
        \tensorho{5}{1};\tensorho{3}{3};\tensorho{3}{4};
        \tensorTD{1}{4};
        \tensorTD{1}{3};\tensorTD{2}{3};
        \tensorTD{1}{1};\tensorTD{2}{1};\tensorTD{4}{1};
        \draw[|<->|](0.5,5.5)--(4.5,5.5);\node at(2.5,5.5)[above]{$m$ columns};
        \draw[|<->|](5.5,0.5)--(5.5,4.5);\node at(5.5,2.5)[right]{$n$ rows};
    \end{tikzpicture}
}

\newcommand{\Ti}{
\raisebox{-0.25cm}{
    \begin{tikzpicture}
        \draw(-0.5,-0.5)grid(3.5,0.5);
        \tensorT{0}{0};\tensorT{1}{0};\tensorT{2}{0};\tensorT{3}{0};
        \node at(0,0.5)[above]{$a_1$};\node at(1,0.5)[above]{$a_2$};\node at(2,0.5)[above]{$a_3$};\node at(3,0.5)[above]{$a_4$};
    \end{tikzpicture}}
}

\newcommand{\Tii}{
\raisebox{-0.25cm}{
    \begin{tikzpicture}
        \draw(-0.5,-0.5)grid(1.5,0.5);\draw(0.33,-0.5)--(0.33,0.5);\draw(0.67,-0.5)--(0.67,0.5);
        \filldraw[fill=gray!20](-0.25,-0.25)rectangle(1.25,0.25);
        \node at(0.5,0) {$\varT^{(4)}$};
        \node at(0.5,0.5)[above]{$a_1 a_2 a_3 a_4$};
    \end{tikzpicture}}
}

\newcommand{\Tiii}{
\raisebox{-0.25cm}{
    \begin{tikzpicture}
        \draw(-0.5,0)--(1.5,0);\draw(0.5,-0.5)--(0.5,0.5);
        \filldraw[fill=gray!20](-0.25,-0.25)rectangle(1.25,0.25);
        \node at(0.5,0) {$\varT^{(4)}$};
        \node at(0.5,0.5)[above]{$a$};
    \end{tikzpicture}}
}

\title{Independent Set Enumeration and Estimation of Related Constants of Grid Graphs and Their Variants}
\author{Kai Liang}
\date{\today} 

\begin{document}
    \maketitle

    \begin{abstract}
        We applied tensor network contraction algorithms to compute the hard-core lattice gas model, i.e., the enumeration of independent sets on grid graphs.
        We observed the influence of surface effect and parity effect on the enumeration (and entropy), and derived upper and lower bounds for both the combinatorics entropy and the coefficients of surface effect by numerical analysis.

        Additionally, we conducted corresponding calculations and analyses for triangular grid graphs, king graph, and cylindrical grid graph.
        We computed and analyzed their associated constants and compared how different adjacency and boundary conditions affect these constants.

        Our computational results have contributed substantial new terms to the OEIS sequence A089980, A027740, A219741, A226444, A245013 and A286513. 
        In addition, we have provided fairly accurate estimates of the relevant constants through numerical analysis of the obtained results.
        Among them, our valuation of the hard square entropy constant is more accurate than existing results.
        And we conject that the surface effect of the periodic boundary of the cylindrical grid graph is $0$—its estimated value of coefficients is very close to $0$.

    \end{abstract}
    
    \section{Introduction}

    Gaunt and Fisher~\cite{LatGas1}, and Runnels and Combs~\cite{LatGas2} independently proposed the \textit{hard-core lattice gas model} in 1965 and 1966, respectively: assume that ``particles'' are placed in the squares of an $m\times n$ lattice, and due to steric effect between particles, no two adjacent (horizontally or vertically) positions can both contain particles.
    Clearly, a valid configuration of $c$ particles corresponds one-to-one with an independent (vertex) set of size $c$ (i.e., a vertex subset with $c$ vertices where no two vertices are adjacent) in the $m\times n$ grid graph.
    The following figure is an example for an independent set with $9$ vertices (the \textbf{black} vertices) of an in the $6\times4$ grid graph:

    \begin{center}
        \squgridi
    \end{center}

    Baxter~\cite{grid1} analyzed this model in 1979 and proved that as the length and width $m$, $n$ of the grid graph tend to infinity, the enumeration of independent sets approaches $(e^f)^{mn}$, where $e^f$ is a constant called the \textit{entropy constant} of the model,
    and its natural logarithm $f$ is the \textit{combinatorial entropy}\footnote{These two constants have many different names in different literature.
    For example, \textit{average combinatorics entropy} or \textit{entropy per unit}, etc.} in the infinite grid graph case, which can be regarded as the total entropy is the production of the entropy distributed by each unit, and their product is the total entropy of the system.
    Using the transfer matrix method, Calkin and Wilf~\cite{grid2} calculated that:
    \begin{equation}
        1.50304808 \leq e^f\leq 1.51316067
    \end{equation}

    Subsequent studies~\cite{grid3, CTM2, grid5, CTM, grid4} provided more precise estimates of this constant:
    \begin{equation}
        e^f=1.5030480824753322643220663294755536893857810..\ldots
    \end{equation}

    Additionally, the lattice gas model has several simple variants.
    For example, replacing the $m\times n$ grid graph with a triangular grid graph of width $n$ (equivalent to substituting the square lattices with regular hexagons).
    The following figure is an example for an independent set with $6$ vertices in $6$-triangular grid graph (i.e., the triangular grid graph with three sides of length $6$):

    \begin{center}
        \hextrigridi
    \end{center}

    Interestingly, the entropy constant for this variant was proven to be an algebraic number by using number theoretic methods,
    specifically, the unique real root of an irreducible 24th-degree polynomial~\cite{hexaConst,hexa1,hexa2}.
    However, no similar results exist for the (original) grid graph or its other conventional variants.
    Only estimated values or upper and lower bounds can be provided.

    We can also change the boundary type of the grid graph to get its variation.
    For example, if the horizontal \textit{closed boundary} is changed to \textit{periodic boundary}, that is, the vertices of the left and right sides are also connected according to the grid pattern, the \textit{cylindrical grid graph} can be obtained.
    The following figures are examples of the independent sets with $10$ vertices on a $6\times 4$ cylindrical grid graph:
    \begin{equation}
        \cylgridi
    \end{equation}

    Berger~\cite{WT_und_01} has demonstrated that no efficient algorithm exists for exact computation in large $m$ and $n$ cases of similar problems.
    This makes the approximation for large-scale cases necessary, which certainly relies on the estimation of relevant constants such as entropy constant.

    Calkin and Wilf~\cite{grid2} established inequalities relating the constants derived from finite configurations under two boundary conditions to the entropy constant in the infinite case, which can be used to provide the strict upper and lower bounds for the entropy constant.
    They provided fairly accurate upper and lower bounds for the hard-square entropy constant by calculating the eigenvalues of the transition matrix for smaller width cases (including cylindrical cases).
    Subsequently, Lundow, Markstr\"om and Friedland~\cite{grid5, grid6} optimized the calculation of the method, and calculated more accurate upper and lower bounds.
    However, this approach requires computing the principal eigenvalue of the transfer matrix, whose dimension grows exponentially with the width $m$.

    An another general method for estimating the entropy constants of such models involves expressing the independent set enumeration as a two-dimensional tensor network contraction and computing or approximating it using the corner transfer matrix method (CTM) and its variants~\cite{grid1, CTM2, grid4}.
    Chan and Rechnitze~\cite{lower, upper} gave the lower and upper bounds of the entropy constants using the CTM renormalization group (CTMRG) method.

    In addition, Ronnie~\cite{prob} used more complex probabilistic methods to estimate the entropy constant, and proved that it is possible to approximate the constant with an additive error of $\epsilon$ with computational time $\left(\frac{1}{\epsilon}\right)^{O(1)}$.

    In this paper, we provide exact enumeration of independent sets for small-scale (but the range is greater than the existing results) grid graphs via two-dimensional tensor network contraction.
    We also calculate the triangular grid graphs, king graphs and cylindrical grid graph.
    Our computational results extend multiple OEIS sequences defined as independent set enumeration as shown in Table~\ref{tab:oeis}.
    \begin{table}
        \small
        \centering
        \begin{tabular}{|c|c|c|c|}
            \hline
            sequence & graph type  & original range & new range \\
            \hline
            A089934 &                    & $1\ldots1225$ & $1\ldots3003$ \\
            A089980 & grid graph         & $(m+n\leq50)$ & $(m+n\leq78)$ \\
            \hline
                    & $n$-triangular     & $n=0\ldots17$ & $n=0\ldots40$ \\
            A027740 & grid graph         & $(n\leq17)$   & $(n\leq40)$   \\
            \hline
            A226444 & $(m,n)$-triangular & $1\ldots1034$ & $1\ldots2926$ \\
            A219741 & grid graph         & $(m+n\leq44)$ & $(m+n\leq77)$ \\
            \hline
                    &                    & $1\ldots1081$ & $1\ldots3320$ \\
            A245013 & king graph         & $(m+n\leq43)$ & $(m+n\leq78)$ \\
            \hline
                    & cylindrical        & $1\ldots435$  & $1\ldots704$  \\
            A286513 & grid graph         & $(m+n\leq29)$ & $(m+n\leq37)$ \\
            \hline
        \end{tabular}
        \normalsize
        \caption{\small Related OEIS sequences. Some of the sequences contain essentially the same data.}
        \label{tab:oeis}
    \end{table}
    We observed the results of these finite cases and found that the difference in entropy per unit compared to the infinite case is mainly influenced by two factors: the parity of the boundary length and the density of the boundary, namely the parity effect and the surface effect.
    When both the length and width are large enough, the surface effect will dominate, which is approximately proportional to the boundary density.
    Therefore, we can estimate the coefficients of surface effect with respect to boundary density based on existing data, and use this to provide an estimate of unit area entropy for infinite cases of each model (as shown in Table~\ref{tab:f}).
    These two constants can be used to provide a more accurate approximation of the unit area entropy or total entropy for finite cases.

    Although our estimates of these constants rely on empirical observations of patterns in small-scale cases, some of the results obtained (such as the hard-square entropy constant) prove to be more accurate than existing ones.

    \begin{table}
        \small
        \centering
        \begin{tabular}{|l|}
            \hline
            grid graph:  \\
            \hspace{2em}$\underline{0.407495101260688000450146812~9046} < f<
            \underline{0.407495101260688000450146812~3584};$ \\
            \hspace{2em}$\underline{1.5030480824753322643220663~2947} < e^{f} <
            \underline{1.5030480824753322643220663~3030};$ \\
            \hspace{2em}$0.0670552316597760242912072~6308 < k<
            0.0670552316597760242912072~7191.$\\
            \hline
            $(m,n)$-triangular grid graph:\\
            \hspace{2em}$\underline{0.33324272197618188785374776~3953} <f_\varDelta <
            \underline{0.33324272197618188785374776~4006};$\\
            \hspace{2em}$\underline{1.395485972479302735229500663~4939} <e^{f_\varDelta}<
            \underline{1.395485972479302735229500663~5675};$\\
            \hspace{2em}$0.11182308857658~5958 <k_\varDelta <
            0.11182308857658~6865.$\\
            \hline
            king graph:\\
            \hspace{2em}$\underline{0.294640767816144}918~4083 <f_K <
            \underline{0.294640767816144918~2294};$\\
            \hspace{2em}$\underline{1.34264395112460}129~7851 <e^{f_K}<
            \underline{1.34264395112460129~8092};$\\
            \hspace{2em}$0.1354918026~4626 <k_K <
            0.1354918026~7005.$\\
            \hline
            cylindrical grid graph:\\
            \hspace{2em}$f_C =f; e^{f_C}=e^f;\bar{k}_C=k;$ \\
            \hspace{2em}$-0.00000000000~1581 < \mathring{k}_C <
             0.00000000000~3415;$ \\
            \hspace{2em}$(\text{conjecture:} ~~\mathring{k}_C=0)$\\
            \hline
            3-dimensional grid graph:\\
            \hspace{2em}$0.3622~2260 < f_{\text{3d}} < 0.3622~8469;$ \\
            \hspace{2em}$1.436~5186 < e^{f_{\text{3d}}} < 1.436~6079;$ \\
            \hspace{2em}$0.04~4592 < k_{\text{3d}} < 0.04~5172.$\\
            \hline
        \end{tabular}
        \normalsize
        \caption{\small Results of constant estimates. The \underline{underlined} parts are the existing results.}
        \label{tab:f}
    \end{table}
%

    \section{Grid graph}

    We denote the number of independent sets in an $m \times n$ grid graph as $N(m,n)$.

    \subsection{Algorithms}

    In both this and the next section, we further optimized the algorithm used in Reference~\cite{my1}, expanding the range of computable results from $\min(m, n)\leq 27$ to $\min(m,n)\leq 39$.

    Baxter~\cite{grid1} provided a method to transform the enumeration of independent sets in grid graphs into a two-dimensional tensor network contraction.
    For convenience, we adopt an equivalent formulation here:
    select (\textbf{bold}) the edge above and to the right of each selected point (add ``half an edge'' if there is no edge), then all consecutive vertices along the same diagonal will be connected, forming non-overlapping W-shaped zigzag paths.
    It is easy to observe that such a set of paths in an $m \times n$ grid graph corresponds bijectively to an independent vertex set.

    \begin{center}
        \gridi~~\raisebox{0.9cm}{$\Longrightarrow$}~~\gridii
    \end{center}

    Furthermore, the edge set in the graph forms such a set of paths if and only if the connections around each vertex follow one of the following five templates
    (note that these templates can only assemble into W-shaped zigzag paths starting at $\wtilei{}{}{}{}$ and ending at $\wtileiv{}{}{}{}$):
    
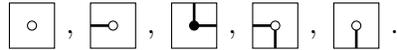
\vspace{-1em}
    \[
        \wtileo{}{}{}{}~,~\wtilei{}{}{}{}~,~\wtileii{}{}{}{}~,~\wtileiii{}{}{}{}~,~\wtileiv{}{}{}{}~.
    \]
    Therefore, we define the index set $\{0, 1\}$, where the indices $0$ and $1$ represent the absence or presence of a connecting edge, respectively.
    The tensor $\varT$ of shape $2\times2\times2\times2$ is constructed based on these five templates:
    \begin{equation}\label{eq:T}
        \varT\coloneq\wtileo{0}{0}{0}{0}+\wtilei{1}{0}{0}{0}+\wtileii{0}{1}{0}{1}+\wtileiii{1}{0}{1}{0}+\wtileiv{0}{0}{1}{0}.
    \end{equation}

    Thus, $N(m,n)$ can be expressed as the contraction of the following two-dimensional tensor network:
    \begin{equation}\label{eq:TN}
        N(m,n)=\raisebox{-3cm}{\TN}
    \end{equation}

    Note that the tensors on the upper and right boundaries contain two indices: 0 and 1, which correspond to the two cases of no path and extended path respectively.

    To reduce the memory usage, we first merge $l$ adjacent tensors horizontally (we call this $l$\textit{-merging}, and $l$ the \textit{stride}), with residual tensors merged if there are on the right most end.
    In this paper we always take $l = 4$, thus,
    \begin{equation}
        \Ti \raisebox{0.25cm}{=} \Tii \raisebox{0.25cm}{=} \Tiii.
    \end{equation}
    After this contraction, the vertical indices $a_1 a_2 a_3 a_4$ become sequences of characters 0s and 1s of length $4$.
    However, only 8 (which is the size of index $a$) out of the 16 possible sequences actually appear, since no two adjacent entries in a sequence can both be 1 (otherwise, the corresponding paths would overlap).
    This significantly reduces the size of the state tensors that need to be stored during row-by-row contraction.
    
    The space efficiency of this merging approach stems from the inherent sparsity of the state tensors in our algorithm. 
    For a given width $m$, each term in the computational state tensor corresponds to a character string of length $m+1$, comprising $m$ vertical characters and a horizontal character (which may appear at any position).
    After performing tensor contraction on each row, the corresponding term become zero when any two adjacent vertical characters are both 1.
    The maximum number of non-zero terms is
    \begin{equation}
        M\coloneq F_{n+2},
    \end{equation}
    the $n+2$-th Fibonacci number.
    
    And after $l$-merging, the maximum size (including zero terms) of state tensors during the contraction process becomes:
    \begin{equation}
        M(l)\coloneq
        \begin{cases}
            F_{n+1}, &m\leq l;\\
            2\cdot F_{l+2}{}^q \cdot F_{r+2}, & m>l,
        \end{cases}
        \text{~where~}n=ql+r, 0<r\leq l.
    \end{equation}
    Here, $F_{l+2}$ represents the number of remained vertical indices after $l$-merging, which equivalent to binary sequences of length $l$ without consecutive 1s;
    $F_{r+2}$ is the number of remained vertical indices of rightmost tensors.
    
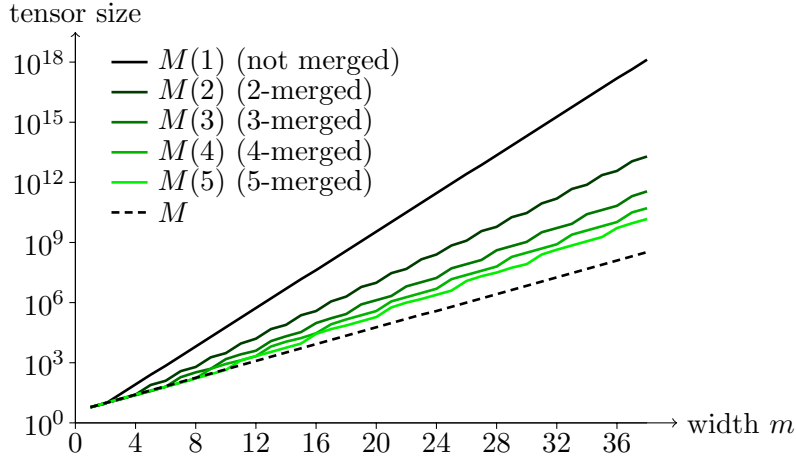
\begin{figure}
         \centering

         \begin{tikzpicture}[scale=0.8]
            \draw[->](0,0)--(10,0);\draw(0,-0.1)grid(9.5,0);\node at(0,0)[below]{0};\node at(1,0)[below]{4};\node at(2,0)[below]{8};\node at(3,0)[below]{12};\node at(4,0)[below]{16};\node at(5,0)[below]{20};\node at(6,0)[below]{24};\node at(7,0)[below]{28};\node at(8,0)[below]{32};\node at(9,0)[below]{36};\node at(10,0)[right]{width $m$}; 
            \draw[->](0,0)--(0,6.5);\draw(-0.1,0)grid(0,6.5);\node at(0,0)[left]{$10^{0}$};\node at(0,1)[left]{$10^{3}$};\node at(0,2)[left]{$10^{6}$};\node at(0,3)[left]{$10^{9}$};\node at(0,4)[left]{$10^{12}$};\node at(0,5)[left]{$10^{15}$};\node at(0,6)[left]{$10^{18}$};\node at(0,6.5)[above]{tensor size}; 

            \draw[line width=1pt](0.25,0.26)--(0.5,0.32)--(0.75,0.48)--(1.0,0.64)--(1.25,0.8)--(1.5,0.95)--(1.75,1.11)--(2.0,1.27)--(2.25,1.43)--(2.5,1.59)--(2.75,1.75)--(3.0,1.91)--(3.25,2.07)--(3.5,2.23)--(3.75,2.39)--(4.0,2.54)--(4.25,2.7)--(4.5,2.86)--(4.75,3.02)--(5.0,3.18)--(5.25,3.34)--(5.5,3.5)--(5.75,3.66)--(6.0,3.82)--(6.25,3.98)--(6.5,4.14)--(6.75,4.29)--(7.0,4.45)--(7.25,4.61)--(7.5,4.77)--(7.75,4.93)--(8.0,5.09)--(8.25,5.25)--(8.5,5.41)--(8.75,5.57)--(9.0,5.73)--(9.25,5.88)--(9.5,6.04); 
            \draw[line width=1pt,color=colori](0.25,0.26)--(0.5,0.33)--(0.75,0.39)--(1.0,0.47)--(1.25,0.63)--(1.5,0.7)--(1.75,0.86)--(2.0,0.93)--(2.25,1.09)--(2.5,1.16)--(2.75,1.32)--(3.0,1.4)--(3.25,1.56)--(3.5,1.63)--(3.75,1.79)--(4.0,1.86)--(4.25,2.02)--(4.5,2.1)--(4.75,2.26)--(5.0,2.33)--(5.25,2.49)--(5.5,2.56)--(5.75,2.72)--(6.0,2.8)--(6.25,2.95)--(6.5,3.03)--(6.75,3.19)--(7.0,3.26)--(7.25,3.42)--(7.5,3.49)--(7.75,3.65)--(8.0,3.73)--(8.25,3.89)--(8.5,3.96)--(8.75,4.12)--(9.0,4.19)--(9.25,4.35)--(9.5,4.43); 
            \draw[line width=1pt,color=colorii](0.25,0.26)--(0.5,0.33)--(0.75,0.4)--(1.0,0.46)--(1.25,0.53)--(1.5,0.6)--(1.75,0.76)--(2.0,0.84)--(2.25,0.9)--(2.5,1.06)--(2.75,1.14)--(3.0,1.2)--(3.25,1.36)--(3.5,1.44)--(3.75,1.51)--(4.0,1.66)--(4.25,1.74)--(4.5,1.81)--(4.75,1.97)--(5.0,2.04)--(5.25,2.11)--(5.5,2.27)--(5.75,2.34)--(6.0,2.41)--(6.25,2.57)--(6.5,2.64)--(6.75,2.71)--(7.0,2.87)--(7.25,2.94)--(7.5,3.01)--(7.75,3.17)--(8.0,3.24)--(8.25,3.31)--(8.5,3.47)--(8.75,3.54)--(9.0,3.61)--(9.25,3.77)--(9.5,3.85); 
            \draw[line width=1pt,color=coloriii](0.25,0.26)--(0.5,0.33)--(0.75,0.4)--(1.0,0.47)--(1.25,0.53)--(1.5,0.6)--(1.75,0.67)--(2.0,0.74)--(2.25,0.9)--(2.5,0.98)--(2.75,1.04)--(3.0,1.11)--(3.25,1.27)--(3.5,1.35)--(3.75,1.41)--(4.0,1.49)--(4.25,1.64)--(4.5,1.72)--(4.75,1.79)--(5.0,1.86)--(5.25,2.02)--(5.5,2.09)--(5.75,2.16)--(6.0,2.23)--(6.25,2.39)--(6.5,2.46)--(6.75,2.53)--(7.0,2.6)--(7.25,2.76)--(7.5,2.83)--(7.75,2.9)--(8.0,2.97)--(8.25,3.13)--(8.5,3.2)--(8.75,3.27)--(9.0,3.34)--(9.25,3.5)--(9.5,3.57); 
            \draw[line width=1pt,color=coloriv](0.25,0.26)--(0.5,0.33)--(0.75,0.4)--(1.0,0.47)--(1.25,0.54)--(1.5,0.6)--(1.75,0.67)--(2.0,0.74)--(2.25,0.81)--(2.5,0.88)--(2.75,1.04)--(3.0,1.11)--(3.25,1.18)--(3.5,1.25)--(3.75,1.32)--(4.0,1.48)--(4.25,1.56)--(4.5,1.62)--(4.75,1.69)--(5.0,1.76)--(5.25,1.92)--(5.5,2.0)--(5.75,2.06)--(6.0,2.13)--(6.25,2.2)--(6.5,2.36)--(6.75,2.44)--(7.0,2.5)--(7.25,2.58)--(7.5,2.64)--(7.75,2.8)--(8.0,2.88)--(8.25,2.95)--(8.5,3.02)--(8.75,3.09)--(9.0,3.24)--(9.25,3.32)--(9.5,3.39); 
            \draw[line width=1pt,densely dashed](0.25,0.26)--(0.5,0.33)--(0.75,0.4)--(1.0,0.47)--(1.25,0.54)--(1.5,0.61)--(1.75,0.68)--(2.0,0.75)--(2.25,0.82)--(2.5,0.89)--(2.75,0.96)--(3.0,1.03)--(3.25,1.1)--(3.5,1.17)--(3.75,1.24)--(4.0,1.31)--(4.25,1.38)--(4.5,1.45)--(4.75,1.52)--(5.0,1.59)--(5.25,1.66)--(5.5,1.73)--(5.75,1.8)--(6.0,1.86)--(6.25,1.93)--(6.5,2.0)--(6.75,2.07)--(7.0,2.14)--(7.25,2.21)--(7.5,2.28)--(7.75,2.35)--(8.0,2.42)--(8.25,2.49)--(8.5,2.56)--(8.75,2.63)--(9.0,2.7)--(9.25,2.77)--(9.5,2.84); 

            \draw[line width=1pt](0.6,6.0)--(1.2,6.0);\node at(1.2,6.0)[right]{$M(1)$ (not merged)}; 
            \draw[line width=1pt,color=colori](0.6,5.5)--(1.2,5.5);\node at(1.2,5.5)[right]{$M(2)$ (2-merged)}; 
            \draw[line width=1pt,color=colorii](0.6,5.0)--(1.2,5.0);\node at(1.2,5.0)[right]{$M(3)$ (3-merged)}; 
            \draw[line width=1pt,color=coloriii](0.6,4.5)--(1.2,4.5);\node at(1.2,4.5)[right]{$M(4)$ (4-merged)}; 
            \draw[line width=1pt,color=coloriv](0.6,4.0)--(1.2,4.0);\node at(1.2,4.0)[right]{$M(5)$ (5-merged)}; 
            \draw[line width=1pt,densely dashed](0.6,3.5)--(1.2,3.5);\node at(1.2,3.5)[right]{$M$}; 
        \end{tikzpicture}
         \caption{\small State tensor sizes after $l$-merged. $M\coloneq F_{n+2}$ is the limit, i.e., the entire rows merged.}
         \label{fig:charnum}
    \end{figure}

    Figure~\ref{fig:charnum} demonstrates the logarithmic value of $M$ and $M(l), l=1,2,\ldots5$ for different widths.
    
    It can be observed that larger $l$ reduces state tensor size (approaching $M$), but with diminishing returns.
    In addition, increased $l$ will raise the order of merged tensors, potentially increase computational complexity.
    Therefore, the stride $l$ must be properly selected to achieve maximum efficiency.

    In addition, we can also use sparse matrices to represent the state tensors, as they actually only participate matrix operations in the algorithm.
    However, its actual effect is not as good as $4$-merging, because when the width is not large, such merging is enough to make the state tensor dense.
    For $m=34$ (the maximum width we calculated) as example, the tensor's density increase from $0.01\%$ to $59.33\%$ after $4$-merges.

    On the basis of the above, and using the congruence enumeration to avoid numerical overflow
    , 
    we contracted the above tensor network and computed $N(m,n)$ for all cases where $m+n\leq78$ (added the result of $N(39,39)$ from A006506).
    The complete results have been archived in the corresponding OEIS sequence A089934.
    In addition, we also calculated all cases where $m\leq34, n\leq 100$, for the next constant estimation.

    \subsection{Observation of results}

    We define the \textit{combinatorial entropy} of the finite grid graph cases as:
    \begin{align}\label{eq:f_def}
        f(m,n)   &\coloneq\frac{\ln N(m,n)}{mn}, \\
        f^-(m,n) &\coloneq\frac{\ln N(m,n)}{(m+1)(n+1)}.
    \end{align}
    Here $mn$ and $(m+1)(n+1)$ are the \textit{areas} of the graph under two different definitions.
    It can be considered that the latter includes the area of the boundary (where the total width and height are both $1$), while the former does not.

    We need to consider the area of the boundary because note that there are two ways to assemble $rs$ $m \times n$ grid graphs into a larger grid graph:

    (1) The graphs of $r$ rows and $s$ columns are directly spliced together.
    This can obtain all the independent sets of a $rm \times sn$ grid graph, but the selected vertex set after splicing may not form an independent set.

    (2) Add rows and columns of unselected vertices at all joints to form a $(rm+m-1) \times (sn+n-1)$ grid graph.
    This can ensure that the selected vertices can form an independent set after splicing, but it can not get all the independent sets of the a $(rm+m-1) \times (sn+n-1)$ grid graph.

    So when $m,n\rightarrow\infty$, $f^-(m,n)$ and $f(m,n)$ approach to the constant $f$ from the lower and upper sides respectively:
    \begin{equation}\label{eq:ineq}
        f^-(m,n) < f < f(m,n);
    \end{equation}

    \begin{equation}
        \lim_{m,n\rightarrow\infty}f^-(m,n) =f = \lim_{m,n\rightarrow\infty}f(m,n).
    \end{equation}
    Which can be used to give the lower and upper bounds of $f$.
    Taking $m=34$ and $n=100$, based on our results\footnote{In this article, we always provide four non-exact digits for the resulted upper and lower bounds (after the spaces)}:
    \begin{equation}
        0.2~8852 < f < 0.2~9999.
    \end{equation}
    It can be observed that the upper and lower bounds given in this way are quite imprecise.
    Therefore, we will observe the numerical patterns of these small-scale cases to obtain a more precise estimation.

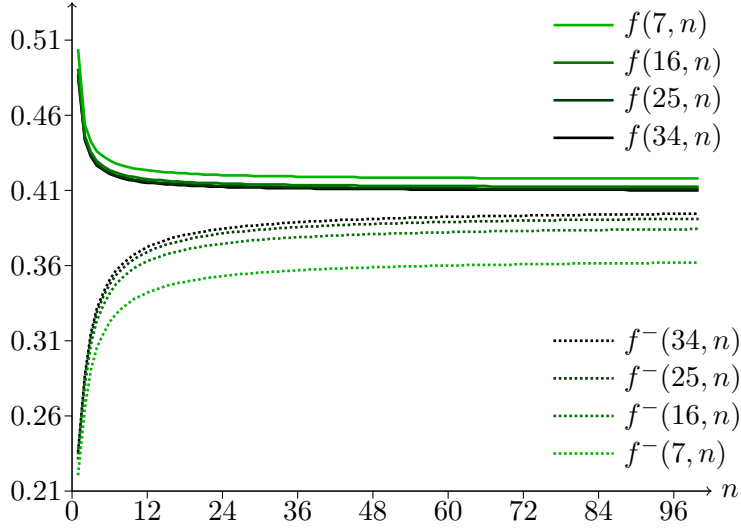
\begin{figure}
         \centering
         \begin{tikzpicture}[scale=1]
            \draw[->](0,0)--(8.5,0);\draw(0,-0.1)grid(8.5,0);\node at(0,0)[below]{0};\node at(1,0)[below]{12};\node at(2,0)[below]{24};\node at(3,0)[below]{36};\node at(4,0)[below]{48};\node at(5,0)[below]{60};\node at(6,0)[below]{72};\node at(7,0)[below]{84};\node at(8,0)[below]{96};\node at(8.5,0)[right]{$n$}; 
            \draw[->](0,0)--(0,6.5);\draw(-0.1,0)grid(0,6.5);\node at(0,0)[left]{0.21};\node at(0,1)[left]{0.26};\node at(0,2)[left]{0.31};\node at(0,3)[left]{0.36};\node at(0,4)[left]{0.41};\node at(0,5)[left]{0.46};\node at(0,6)[left]{0.51};\node at(0, 6.5)[left]{}; 
            \draw[line width=1pt](0.08,5.52)--(0.17,4.67)--(0.25,4.45)--(0.33,4.33)--(0.42,4.27)--(0.5,4.22)--(0.58,4.19)--(0.67,4.16)--(0.75,4.14)--(0.83,4.13)--(0.92,4.11)--(1.0,4.1)--(1.08,4.1)--(1.17,4.09)--(1.25,4.08)--(1.33,4.08)--(1.42,4.07)--(1.5,4.07)--(1.58,4.06)--(1.67,4.06)--(1.75,4.06)--(1.83,4.05)--(1.92,4.05)--(2.0,4.05)--(2.08,4.04)--(2.17,4.04)--(2.25,4.04)--(2.33,4.04)--(2.42,4.04)--(2.5,4.04)--(2.58,4.03)--(2.67,4.03)--(2.75,4.03)--(2.83,4.03)--(2.92,4.03)--(3.0,4.03)--(3.08,4.03)--(3.17,4.03)--(3.25,4.02)--(3.33,4.02)--(3.42,4.02)--(3.5,4.02)--(3.58,4.02)--(3.67,4.02)--(3.75,4.02)--(3.83,4.02)--(3.92,4.02)--(4.0,4.02)--(4.08,4.02)--(4.17,4.02)--(4.25,4.02)--(4.33,4.02)--(4.42,4.02)--(4.5,4.01)--(4.58,4.01)--(4.67,4.01)--(4.75,4.01)--(4.83,4.01)--(4.92,4.01)--(5.0,4.01)--(5.08,4.01)--(5.17,4.01)--(5.25,4.01)--(5.33,4.01)--(5.42,4.01)--(5.5,4.01)--(5.58,4.01)--(5.67,4.01)--(5.75,4.01)--(5.83,4.01)--(5.92,4.01)--(6.0,4.01)--(6.08,4.01)--(6.17,4.01)--(6.25,4.01)--(6.33,4.01)--(6.42,4.01)--(6.5,4.01)--(6.58,4.01)--(6.67,4.01)--(6.75,4.01)--(6.83,4.01)--(6.92,4.01)--(7.0,4.01)--(7.08,4.01)--(7.17,4.01)--(7.25,4.01)--(7.33,4.01)--(7.42,4.0)--(7.5,4.0)--(7.58,4.0)--(7.67,4.0)--(7.75,4.0)--(7.83,4.0)--(7.92,4.0)--(8.0,4.0)--(8.08,4.0)--(8.17,4.0)--(8.25,4.0)--(8.33,4.0); 
            \draw[line width=1pt](6.4,4.7)--(7.2,4.7);\node at(7.2,4.7)[right]{$f(34,n)$}; 
            \draw[line width=1pt,color=colori](0.08,5.55)--(0.17,4.69)--(0.25,4.47)--(0.33,4.35)--(0.42,4.28)--(0.5,4.24)--(0.58,4.2)--(0.67,4.18)--(0.75,4.16)--(0.83,4.14)--(0.92,4.13)--(1.0,4.12)--(1.08,4.11)--(1.17,4.1)--(1.25,4.1)--(1.33,4.09)--(1.42,4.09)--(1.5,4.08)--(1.58,4.08)--(1.67,4.07)--(1.75,4.07)--(1.83,4.07)--(1.92,4.06)--(2.0,4.06)--(2.08,4.06)--(2.17,4.06)--(2.25,4.06)--(2.33,4.05)--(2.42,4.05)--(2.5,4.05)--(2.58,4.05)--(2.67,4.05)--(2.75,4.05)--(2.83,4.04)--(2.92,4.04)--(3.0,4.04)--(3.08,4.04)--(3.17,4.04)--(3.25,4.04)--(3.33,4.04)--(3.42,4.04)--(3.5,4.04)--(3.58,4.04)--(3.67,4.04)--(3.75,4.03)--(3.83,4.03)--(3.92,4.03)--(4.0,4.03)--(4.08,4.03)--(4.17,4.03)--(4.25,4.03)--(4.33,4.03)--(4.42,4.03)--(4.5,4.03)--(4.58,4.03)--(4.67,4.03)--(4.75,4.03)--(4.83,4.03)--(4.92,4.03)--(5.0,4.03)--(5.08,4.03)--(5.17,4.03)--(5.25,4.03)--(5.33,4.03)--(5.42,4.02)--(5.5,4.02)--(5.58,4.02)--(5.67,4.02)--(5.75,4.02)--(5.83,4.02)--(5.92,4.02)--(6.0,4.02)--(6.08,4.02)--(6.17,4.02)--(6.25,4.02)--(6.33,4.02)--(6.42,4.02)--(6.5,4.02)--(6.58,4.02)--(6.67,4.02)--(6.75,4.02)--(6.83,4.02)--(6.92,4.02)--(7.0,4.02)--(7.08,4.02)--(7.17,4.02)--(7.25,4.02)--(7.33,4.02)--(7.42,4.02)--(7.5,4.02)--(7.58,4.02)--(7.67,4.02)--(7.75,4.02)--(7.83,4.02)--(7.92,4.02)--(8.0,4.02)--(8.08,4.02)--(8.17,4.02)--(8.25,4.02)--(8.33,4.02); 
            \draw[line width=1pt,color=colori](6.4,5.2)--(7.2,5.2);\node at(7.2,5.2)[right]{$f(25,n)$}; 
            \draw[line width=1pt,color=colorii](0.08,5.62)--(0.17,4.73)--(0.25,4.51)--(0.33,4.39)--(0.42,4.32)--(0.5,4.27)--(0.58,4.24)--(0.67,4.21)--(0.75,4.19)--(0.83,4.18)--(0.92,4.16)--(1.0,4.15)--(1.08,4.14)--(1.17,4.14)--(1.25,4.13)--(1.33,4.12)--(1.42,4.12)--(1.5,4.11)--(1.58,4.11)--(1.67,4.1)--(1.75,4.1)--(1.83,4.1)--(1.92,4.1)--(2.0,4.09)--(2.08,4.09)--(2.17,4.09)--(2.25,4.09)--(2.33,4.08)--(2.42,4.08)--(2.5,4.08)--(2.58,4.08)--(2.67,4.08)--(2.75,4.08)--(2.83,4.08)--(2.92,4.07)--(3.0,4.07)--(3.08,4.07)--(3.17,4.07)--(3.25,4.07)--(3.33,4.07)--(3.42,4.07)--(3.5,4.07)--(3.58,4.07)--(3.67,4.07)--(3.75,4.07)--(3.83,4.06)--(3.92,4.06)--(4.0,4.06)--(4.08,4.06)--(4.17,4.06)--(4.25,4.06)--(4.33,4.06)--(4.42,4.06)--(4.5,4.06)--(4.58,4.06)--(4.67,4.06)--(4.75,4.06)--(4.83,4.06)--(4.92,4.06)--(5.0,4.06)--(5.08,4.06)--(5.17,4.06)--(5.25,4.06)--(5.33,4.06)--(5.42,4.06)--(5.5,4.06)--(5.58,4.05)--(5.67,4.05)--(5.75,4.05)--(5.83,4.05)--(5.92,4.05)--(6.0,4.05)--(6.08,4.05)--(6.17,4.05)--(6.25,4.05)--(6.33,4.05)--(6.42,4.05)--(6.5,4.05)--(6.58,4.05)--(6.67,4.05)--(6.75,4.05)--(6.83,4.05)--(6.92,4.05)--(7.0,4.05)--(7.08,4.05)--(7.17,4.05)--(7.25,4.05)--(7.33,4.05)--(7.42,4.05)--(7.5,4.05)--(7.58,4.05)--(7.67,4.05)--(7.75,4.05)--(7.83,4.05)--(7.92,4.05)--(8.0,4.05)--(8.08,4.05)--(8.17,4.05)--(8.25,4.05)--(8.33,4.05); 
            \draw[line width=1pt,color=colorii](6.4,5.7)--(7.2,5.7);\node at(7.2,5.7)[right]{$f(16,n)$}; 
            \draw[line width=1pt,color=coloriii](0.08,5.88)--(0.17,4.88)--(0.25,4.65)--(0.33,4.52)--(0.42,4.45)--(0.5,4.4)--(0.58,4.36)--(0.67,4.33)--(0.75,4.31)--(0.83,4.29)--(0.92,4.28)--(1.0,4.27)--(1.08,4.26)--(1.17,4.25)--(1.25,4.24)--(1.33,4.24)--(1.42,4.23)--(1.5,4.23)--(1.58,4.22)--(1.67,4.22)--(1.75,4.21)--(1.83,4.21)--(1.92,4.21)--(2.0,4.2)--(2.08,4.2)--(2.17,4.2)--(2.25,4.2)--(2.33,4.2)--(2.42,4.19)--(2.5,4.19)--(2.58,4.19)--(2.67,4.19)--(2.75,4.19)--(2.83,4.19)--(2.92,4.19)--(3.0,4.18)--(3.08,4.18)--(3.17,4.18)--(3.25,4.18)--(3.33,4.18)--(3.42,4.18)--(3.5,4.18)--(3.58,4.18)--(3.67,4.18)--(3.75,4.18)--(3.83,4.17)--(3.92,4.17)--(4.0,4.17)--(4.08,4.17)--(4.17,4.17)--(4.25,4.17)--(4.33,4.17)--(4.42,4.17)--(4.5,4.17)--(4.58,4.17)--(4.67,4.17)--(4.75,4.17)--(4.83,4.17)--(4.92,4.17)--(5.0,4.17)--(5.08,4.17)--(5.17,4.17)--(5.25,4.17)--(5.33,4.17)--(5.42,4.16)--(5.5,4.16)--(5.58,4.16)--(5.67,4.16)--(5.75,4.16)--(5.83,4.16)--(5.92,4.16)--(6.0,4.16)--(6.08,4.16)--(6.17,4.16)--(6.25,4.16)--(6.33,4.16)--(6.42,4.16)--(6.5,4.16)--(6.58,4.16)--(6.67,4.16)--(6.75,4.16)--(6.83,4.16)--(6.92,4.16)--(7.0,4.16)--(7.08,4.16)--(7.17,4.16)--(7.25,4.16)--(7.33,4.16)--(7.42,4.16)--(7.5,4.16)--(7.58,4.16)--(7.67,4.16)--(7.75,4.16)--(7.83,4.16)--(7.92,4.16)--(8.0,4.16)--(8.08,4.16)--(8.17,4.16)--(8.25,4.16)--(8.33,4.16); 
            \draw[line width=1pt,color=coloriii](6.4,6.2)--(7.2,6.2);\node at(7.2,6.2)[right]{$f(7,n)$}; 
            \draw[line width=1pt,densely dotted](0.08,0.52)--(0.17,1.54)--(0.25,2.1)--(0.33,2.43)--(0.42,2.65)--(0.5,2.81)--(0.58,2.93)--(0.67,3.02)--(0.75,3.09)--(0.83,3.15)--(0.92,3.2)--(1.0,3.25)--(1.08,3.28)--(1.17,3.31)--(1.25,3.34)--(1.33,3.37)--(1.42,3.39)--(1.5,3.41)--(1.58,3.42)--(1.67,3.44)--(1.75,3.45)--(1.83,3.47)--(1.92,3.48)--(2.0,3.49)--(2.08,3.5)--(2.17,3.51)--(2.25,3.52)--(2.33,3.53)--(2.42,3.53)--(2.5,3.54)--(2.58,3.55)--(2.67,3.55)--(2.75,3.56)--(2.83,3.57)--(2.92,3.57)--(3.0,3.58)--(3.08,3.58)--(3.17,3.59)--(3.25,3.59)--(3.33,3.59)--(3.42,3.6)--(3.5,3.6)--(3.58,3.61)--(3.67,3.61)--(3.75,3.61)--(3.83,3.61)--(3.92,3.62)--(4.0,3.62)--(4.08,3.62)--(4.17,3.63)--(4.25,3.63)--(4.33,3.63)--(4.42,3.63)--(4.5,3.64)--(4.58,3.64)--(4.67,3.64)--(4.75,3.64)--(4.83,3.64)--(4.92,3.65)--(5.0,3.65)--(5.08,3.65)--(5.17,3.65)--(5.25,3.65)--(5.33,3.65)--(5.42,3.66)--(5.5,3.66)--(5.58,3.66)--(5.67,3.66)--(5.75,3.66)--(5.83,3.66)--(5.92,3.66)--(6.0,3.66)--(6.08,3.67)--(6.17,3.67)--(6.25,3.67)--(6.33,3.67)--(6.42,3.67)--(6.5,3.67)--(6.58,3.67)--(6.67,3.67)--(6.75,3.67)--(6.83,3.68)--(6.92,3.68)--(7.0,3.68)--(7.08,3.68)--(7.17,3.68)--(7.25,3.68)--(7.33,3.68)--(7.42,3.68)--(7.5,3.68)--(7.58,3.68)--(7.67,3.68)--(7.75,3.68)--(7.83,3.69)--(7.92,3.69)--(8.0,3.69)--(8.08,3.69)--(8.17,3.69)--(8.25,3.69)--(8.33,3.69); 
            \draw[line width=1pt,densely dotted](6.4,2.0)--(7.2,2.0);\node at(7.2,2.0)[right]{$f^-(34,n)$}; 
            \draw[line width=1pt,densely dotted,color=colori](0.08,0.49)--(0.17,1.5)--(0.25,2.05)--(0.33,2.38)--(0.42,2.6)--(0.5,2.75)--(0.58,2.87)--(0.67,2.96)--(0.75,3.03)--(0.83,3.09)--(0.92,3.14)--(1.0,3.18)--(1.08,3.22)--(1.17,3.25)--(1.25,3.28)--(1.33,3.3)--(1.42,3.32)--(1.5,3.34)--(1.58,3.36)--(1.67,3.38)--(1.75,3.39)--(1.83,3.4)--(1.92,3.42)--(2.0,3.43)--(2.08,3.44)--(2.17,3.45)--(2.25,3.45)--(2.33,3.46)--(2.42,3.47)--(2.5,3.48)--(2.58,3.48)--(2.67,3.49)--(2.75,3.5)--(2.83,3.5)--(2.92,3.51)--(3.0,3.51)--(3.08,3.52)--(3.17,3.52)--(3.25,3.52)--(3.33,3.53)--(3.42,3.53)--(3.5,3.54)--(3.58,3.54)--(3.67,3.54)--(3.75,3.55)--(3.83,3.55)--(3.92,3.55)--(4.0,3.55)--(4.08,3.56)--(4.17,3.56)--(4.25,3.56)--(4.33,3.56)--(4.42,3.57)--(4.5,3.57)--(4.58,3.57)--(4.67,3.57)--(4.75,3.58)--(4.83,3.58)--(4.92,3.58)--(5.0,3.58)--(5.08,3.58)--(5.17,3.58)--(5.25,3.59)--(5.33,3.59)--(5.42,3.59)--(5.5,3.59)--(5.58,3.59)--(5.67,3.59)--(5.75,3.59)--(5.83,3.6)--(5.92,3.6)--(6.0,3.6)--(6.08,3.6)--(6.17,3.6)--(6.25,3.6)--(6.33,3.6)--(6.42,3.6)--(6.5,3.61)--(6.58,3.61)--(6.67,3.61)--(6.75,3.61)--(6.83,3.61)--(6.92,3.61)--(7.0,3.61)--(7.08,3.61)--(7.17,3.61)--(7.25,3.61)--(7.33,3.61)--(7.42,3.62)--(7.5,3.62)--(7.58,3.62)--(7.67,3.62)--(7.75,3.62)--(7.83,3.62)--(7.92,3.62)--(8.0,3.62)--(8.08,3.62)--(8.17,3.62)--(8.25,3.62)--(8.33,3.62); 
            \draw[line width=1pt,densely dotted,color=colori](6.4,1.5)--(7.2,1.5);\node at(7.2,1.5)[right]{$f^-(25,n)$}; 
            \draw[line width=1pt,densely dotted,color=colorii](0.08,0.42)--(0.17,1.4)--(0.25,1.95)--(0.33,2.27)--(0.42,2.48)--(0.5,2.63)--(0.58,2.75)--(0.67,2.84)--(0.75,2.91)--(0.83,2.97)--(0.92,3.02)--(1.0,3.06)--(1.08,3.09)--(1.17,3.12)--(1.25,3.15)--(1.33,3.17)--(1.42,3.19)--(1.5,3.21)--(1.58,3.23)--(1.67,3.24)--(1.75,3.26)--(1.83,3.27)--(1.92,3.28)--(2.0,3.29)--(2.08,3.3)--(2.17,3.31)--(2.25,3.32)--(2.33,3.33)--(2.42,3.34)--(2.5,3.34)--(2.58,3.35)--(2.67,3.36)--(2.75,3.36)--(2.83,3.37)--(2.92,3.37)--(3.0,3.38)--(3.08,3.38)--(3.17,3.38)--(3.25,3.39)--(3.33,3.39)--(3.42,3.4)--(3.5,3.4)--(3.58,3.4)--(3.67,3.41)--(3.75,3.41)--(3.83,3.41)--(3.92,3.42)--(4.0,3.42)--(4.08,3.42)--(4.17,3.42)--(4.25,3.43)--(4.33,3.43)--(4.42,3.43)--(4.5,3.43)--(4.58,3.43)--(4.67,3.44)--(4.75,3.44)--(4.83,3.44)--(4.92,3.44)--(5.0,3.44)--(5.08,3.45)--(5.17,3.45)--(5.25,3.45)--(5.33,3.45)--(5.42,3.45)--(5.5,3.45)--(5.58,3.46)--(5.67,3.46)--(5.75,3.46)--(5.83,3.46)--(5.92,3.46)--(6.0,3.46)--(6.08,3.46)--(6.17,3.46)--(6.25,3.47)--(6.33,3.47)--(6.42,3.47)--(6.5,3.47)--(6.58,3.47)--(6.67,3.47)--(6.75,3.47)--(6.83,3.47)--(6.92,3.47)--(7.0,3.47)--(7.08,3.47)--(7.17,3.48)--(7.25,3.48)--(7.33,3.48)--(7.42,3.48)--(7.5,3.48)--(7.58,3.48)--(7.67,3.48)--(7.75,3.48)--(7.83,3.48)--(7.92,3.48)--(8.0,3.48)--(8.08,3.48)--(8.17,3.48)--(8.25,3.49)--(8.33,3.49); 
            \draw[line width=1pt,densely dotted,color=colorii](6.4,1.0)--(7.2,1.0);\node at(7.2,1.0)[right]{$f^-(16,n)$}; 
            \draw[line width=1pt,densely dotted,color=coloriii](0.08,0.21)--(0.17,1.1)--(0.25,1.61)--(0.33,1.91)--(0.42,2.1)--(0.5,2.25)--(0.58,2.35)--(0.67,2.44)--(0.75,2.5)--(0.83,2.56)--(0.92,2.6)--(1.0,2.64)--(1.08,2.67)--(1.17,2.7)--(1.25,2.73)--(1.33,2.75)--(1.42,2.77)--(1.5,2.78)--(1.58,2.8)--(1.67,2.81)--(1.75,2.83)--(1.83,2.84)--(1.92,2.85)--(2.0,2.86)--(2.08,2.87)--(2.17,2.88)--(2.25,2.89)--(2.33,2.89)--(2.42,2.9)--(2.5,2.91)--(2.58,2.91)--(2.67,2.92)--(2.75,2.92)--(2.83,2.93)--(2.92,2.93)--(3.0,2.94)--(3.08,2.94)--(3.17,2.95)--(3.25,2.95)--(3.33,2.95)--(3.42,2.96)--(3.5,2.96)--(3.58,2.96)--(3.67,2.97)--(3.75,2.97)--(3.83,2.97)--(3.92,2.97)--(4.0,2.98)--(4.08,2.98)--(4.17,2.98)--(4.25,2.98)--(4.33,2.99)--(4.42,2.99)--(4.5,2.99)--(4.58,2.99)--(4.67,2.99)--(4.75,3.0)--(4.83,3.0)--(4.92,3.0)--(5.0,3.0)--(5.08,3.0)--(5.17,3.0)--(5.25,3.01)--(5.33,3.01)--(5.42,3.01)--(5.5,3.01)--(5.58,3.01)--(5.67,3.01)--(5.75,3.01)--(5.83,3.01)--(5.92,3.02)--(6.0,3.02)--(6.08,3.02)--(6.17,3.02)--(6.25,3.02)--(6.33,3.02)--(6.42,3.02)--(6.5,3.02)--(6.58,3.02)--(6.67,3.03)--(6.75,3.03)--(6.83,3.03)--(6.92,3.03)--(7.0,3.03)--(7.08,3.03)--(7.17,3.03)--(7.25,3.03)--(7.33,3.03)--(7.42,3.03)--(7.5,3.03)--(7.58,3.03)--(7.67,3.03)--(7.75,3.04)--(7.83,3.04)--(7.92,3.04)--(8.0,3.04)--(8.08,3.04)--(8.17,3.04)--(8.25,3.04)--(8.33,3.04); 
            \draw[line width=1pt,densely dotted,color=coloriii](6.4,0.5)--(7.2,0.5);\node at(7.2,0.5)[right]{$f^-(7,n)$}; 
        \end{tikzpicture}
         \caption{\small $f(m,n), f^-(m,n)$ for some fixed $m$}
         \label{fig:f_mn}
    \end{figure}

\begin{figure}
         \centering
         \begin{tikzpicture}[scale=1]
            \draw[->](0,0)--(8.5,0);\draw(0,-0.1)grid(8.5,0);\node at(0,0)[below]{$\frac{1}{\infty}$};\node at(1,0)[below]{$\frac{1}{8}$};\node at(2,0)[below]{$\frac{1}{4}$};\node at(4,0)[below]{$\frac{1}{2}$};\node at(8,0)[below]{1};\node at(8.5,0)[right]{$m$}; 
            \draw[->](0,0)--(0,6.5);\draw(-0.1,0)grid(0,6.5);\node at(0,0)[left]{0.21};\node at(0,1)[left]{0.26};\node at(0,2)[left]{0.31};\node at(0,3)[left]{0.36};\node at(0,4)[left]{0.41};\node at(0,5)[left]{0.46};\node at(0,6)[left]{0.51}; 
            \draw[line width=1pt,color=coloriii](8.0,5.88)--(4.0,4.88)--(2.67,4.65)--(2.0,4.52)--(1.6,4.45)--(1.33,4.4)--(1.14,4.36)--(1.0,4.33)--(0.89,4.31)--(0.8,4.29)--(0.73,4.28)--(0.67,4.27)--(0.62,4.26)--(0.57,4.25)--(0.53,4.24)--(0.5,4.24)--(0.47,4.23)--(0.44,4.23)--(0.42,4.22)--(0.4,4.22)--(0.38,4.21)--(0.36,4.21)--(0.35,4.21)--(0.33,4.2)--(0.32,4.2)--(0.31,4.2)--(0.3,4.2)--(0.29,4.2)--(0.28,4.19)--(0.27,4.19)--(0.26,4.19)--(0.25,4.19)--(0.24,4.19)--(0.24,4.19); 
            \draw[line width=1pt,color=coloriii](7.4,4.7)--(8.2,4.7);\node at(8.2,4.7)[right]{$f(7,n)$}; 
            \draw[line width=1pt,color=colorii](8.0,5.62)--(4.0,4.73)--(2.67,4.51)--(2.0,4.39)--(1.6,4.32)--(1.33,4.27)--(1.14,4.24)--(1.0,4.21)--(0.89,4.19)--(0.8,4.18)--(0.73,4.16)--(0.67,4.15)--(0.62,4.14)--(0.57,4.14)--(0.53,4.13)--(0.5,4.12)--(0.47,4.12)--(0.44,4.11)--(0.42,4.11)--(0.4,4.1)--(0.38,4.1)--(0.36,4.1)--(0.35,4.1)--(0.33,4.09)--(0.32,4.09)--(0.31,4.09)--(0.3,4.09)--(0.29,4.08)--(0.28,4.08)--(0.27,4.08)--(0.26,4.08)--(0.25,4.08)--(0.24,4.08)--(0.24,4.08); 
            \draw[line width=1pt,color=colorii](7.4,4.2)--(8.2,4.2);\node at(8.2,4.2)[right]{$f(16,n)$}; 
            \draw[line width=1pt,color=colori](8.0,5.55)--(4.0,4.69)--(2.67,4.47)--(2.0,4.35)--(1.6,4.28)--(1.33,4.24)--(1.14,4.2)--(1.0,4.18)--(0.89,4.16)--(0.8,4.14)--(0.73,4.13)--(0.67,4.12)--(0.62,4.11)--(0.57,4.1)--(0.53,4.1)--(0.5,4.09)--(0.47,4.09)--(0.44,4.08)--(0.42,4.08)--(0.4,4.07)--(0.38,4.07)--(0.36,4.07)--(0.35,4.06)--(0.33,4.06)--(0.32,4.06)--(0.31,4.06)--(0.3,4.06)--(0.29,4.05)--(0.28,4.05)--(0.27,4.05)--(0.26,4.05)--(0.25,4.05)--(0.24,4.05)--(0.24,4.04); 
            \draw[line width=1pt,color=colori](7.4,3.7)--(8.2,3.7);\node at(8.2,3.7)[right]{$f(25,n)$}; 
            \draw[line width=1pt](8.0,5.52)--(4.0,4.67)--(2.67,4.45)--(2.0,4.33)--(1.6,4.27)--(1.33,4.22)--(1.14,4.19)--(1.0,4.16)--(0.89,4.14)--(0.8,4.13)--(0.73,4.11)--(0.67,4.1)--(0.62,4.1)--(0.57,4.09)--(0.53,4.08)--(0.5,4.08)--(0.47,4.07)--(0.44,4.07)--(0.42,4.06)--(0.4,4.06)--(0.38,4.06)--(0.36,4.05)--(0.35,4.05)--(0.33,4.05)--(0.32,4.04)--(0.31,4.04)--(0.3,4.04)--(0.29,4.04)--(0.28,4.04)--(0.27,4.04)--(0.26,4.03)--(0.25,4.03)--(0.24,4.03)--(0.24,4.03); 
            \draw[line width=1pt](7.4,3.2)--(8.2,3.2);\node at(8.2,3.2)[right]{$f(34,n)$}; 
            \draw[line width=1pt,densely dotted,color=coloriii](8.0,0.21)--(4.0,1.1)--(2.67,1.61)--(2.0,1.91)--(1.6,2.1)--(1.33,2.25)--(1.14,2.35)--(1.0,2.44)--(0.89,2.5)--(0.8,2.56)--(0.73,2.6)--(0.67,2.64)--(0.62,2.67)--(0.57,2.7)--(0.53,2.73)--(0.5,2.75)--(0.47,2.77)--(0.44,2.78)--(0.42,2.8)--(0.4,2.81)--(0.38,2.83)--(0.36,2.84)--(0.35,2.85)--(0.33,2.86)--(0.32,2.87)--(0.31,2.88)--(0.3,2.89)--(0.29,2.89)--(0.28,2.9)--(0.27,2.91)--(0.26,2.91)--(0.25,2.92)--(0.24,2.92)--(0.24,2.93); 
            \draw[line width=1pt,densely dotted,color=coloriii](7.4,1.2)--(8.2,1.2);\node at(8.2,1.2)[right]{$f^-(7,n)$}; 
            \draw[line width=1pt,densely dotted,color=colorii](8.0,0.42)--(4.0,1.4)--(2.67,1.95)--(2.0,2.27)--(1.6,2.48)--(1.33,2.63)--(1.14,2.75)--(1.0,2.84)--(0.89,2.91)--(0.8,2.97)--(0.73,3.02)--(0.67,3.06)--(0.62,3.09)--(0.57,3.12)--(0.53,3.15)--(0.5,3.17)--(0.47,3.19)--(0.44,3.21)--(0.42,3.23)--(0.4,3.24)--(0.38,3.26)--(0.36,3.27)--(0.35,3.28)--(0.33,3.29)--(0.32,3.3)--(0.31,3.31)--(0.3,3.32)--(0.29,3.33)--(0.28,3.34)--(0.27,3.34)--(0.26,3.35)--(0.25,3.36)--(0.24,3.36)--(0.24,3.37); 
            \draw[line width=1pt,densely dotted,color=colorii](7.4,1.7)--(8.2,1.7);\node at(8.2,1.7)[right]{$f^-(16,n)$}; 
            \draw[line width=1pt,densely dotted,color=colori](8.0,0.49)--(4.0,1.5)--(2.67,2.05)--(2.0,2.38)--(1.6,2.6)--(1.33,2.75)--(1.14,2.87)--(1.0,2.96)--(0.89,3.03)--(0.8,3.09)--(0.73,3.14)--(0.67,3.18)--(0.62,3.22)--(0.57,3.25)--(0.53,3.28)--(0.5,3.3)--(0.47,3.32)--(0.44,3.34)--(0.42,3.36)--(0.4,3.38)--(0.38,3.39)--(0.36,3.4)--(0.35,3.42)--(0.33,3.43)--(0.32,3.44)--(0.31,3.45)--(0.3,3.45)--(0.29,3.46)--(0.28,3.47)--(0.27,3.48)--(0.26,3.48)--(0.25,3.49)--(0.24,3.5)--(0.24,3.5); 
            \draw[line width=1pt,densely dotted,color=colori](7.4,2.2)--(8.2,2.2);\node at(8.2,2.2)[right]{$f^-(25,n)$}; 
            \draw[line width=1pt,densely dotted](8.0,0.52)--(4.0,1.54)--(2.67,2.1)--(2.0,2.43)--(1.6,2.65)--(1.33,2.81)--(1.14,2.93)--(1.0,3.02)--(0.89,3.09)--(0.8,3.15)--(0.73,3.2)--(0.67,3.25)--(0.62,3.28)--(0.57,3.31)--(0.53,3.34)--(0.5,3.37)--(0.47,3.39)--(0.44,3.41)--(0.42,3.42)--(0.4,3.44)--(0.38,3.45)--(0.36,3.47)--(0.35,3.48)--(0.33,3.49)--(0.32,3.5)--(0.31,3.51)--(0.3,3.52)--(0.29,3.53)--(0.28,3.53)--(0.27,3.54)--(0.26,3.55)--(0.25,3.55)--(0.24,3.56)--(0.24,3.57); 
            \draw[line width=1pt,densely dotted](7.4,2.7)--(8.2,2.7);\node at(8.2,2.7)[right]{$f^-(34,n)$}; 
        \end{tikzpicture}
         \caption{\small $f(m,n), f^-(m,n)$ for some fixed $m$, where the horizontal axis is $\frac{1}{m}$ (density of left and right boundary)}
         \label{fig:f_mn2}
    \end{figure}
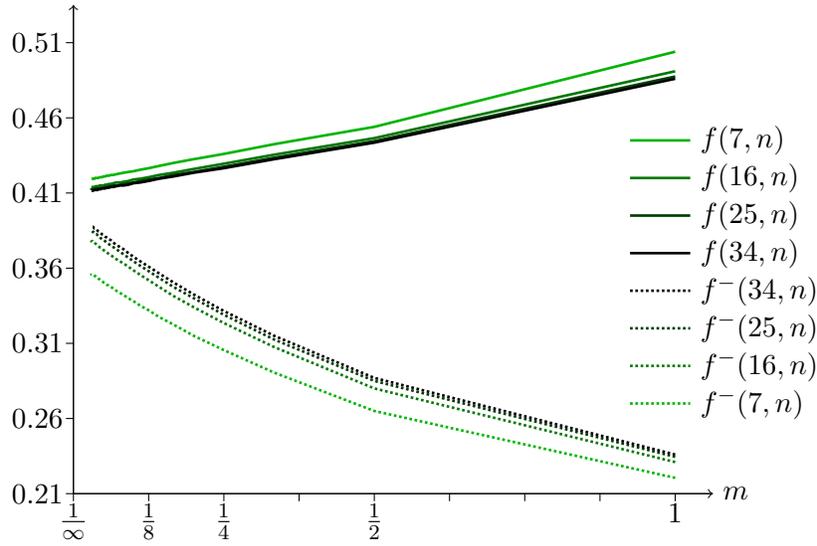

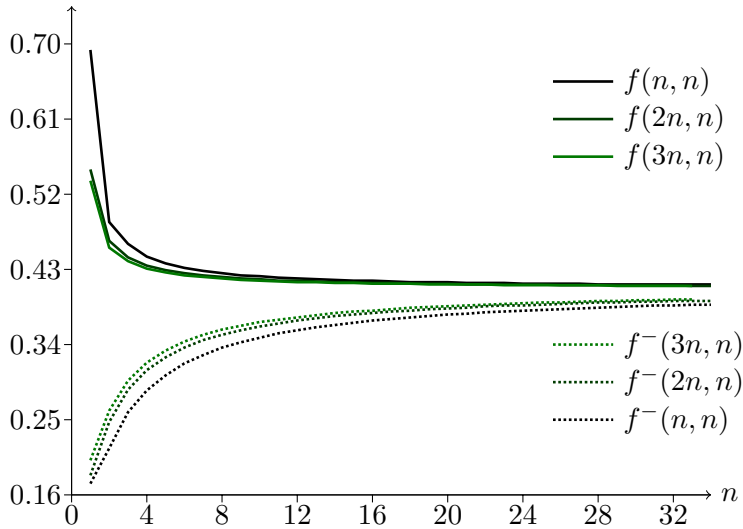
\begin{figure}
         \centering
         \begin{tikzpicture}[scale=1]
            \draw[->](0,0)--(8.5,0);\draw(0,-0.1)grid(8.5,0);\node at(0,0)[below]{0};\node at(1,0)[below]{4};\node at(2,0)[below]{8};\node at(3,0)[below]{12};\node at(4,0)[below]{16};\node at(5,0)[below]{20};\node at(6,0)[below]{24};\node at(7,0)[below]{28};\node at(8,0)[below]{32};\node at(8.5,0)[right]{$n$}; 
            \draw[->](0,0)--(0,6.5);\draw(-0.1,0)grid(0,6.5);\node at(0,0)[left]{0.16};\node at(0,1)[left]{0.25};\node at(0,2)[left]{0.34};\node at(0,3)[left]{0.43};\node at(0,4)[left]{0.52};\node at(0,5)[left]{0.61};\node at(0,6)[left]{0.70};\node at(0, 6.5)[left]{}; 
            \draw[line width=1pt](0.25,5.92)--(0.5,3.63)--(0.75,3.34)--(1.0,3.17)--(1.25,3.08)--(1.5,3.02)--(1.75,2.98)--(2.0,2.95)--(2.25,2.92)--(2.5,2.91)--(2.75,2.89)--(3.0,2.88)--(3.25,2.87)--(3.5,2.86)--(3.75,2.85)--(4.0,2.85)--(4.25,2.84)--(4.5,2.83)--(4.75,2.83)--(5.0,2.83)--(5.25,2.82)--(5.5,2.82)--(5.75,2.82)--(6.0,2.81)--(6.25,2.81)--(6.5,2.81)--(6.75,2.81)--(7.0,2.8)--(7.25,2.8)--(7.5,2.8)--(7.75,2.8)--(8.0,2.8)--(8.25,2.8)--(8.5,2.8); 
            \draw[line width=1pt](6.4,5.5)--(7.2,5.5);\node at(7.2,5.5)[right]{$f(n,n)$}; 
            \draw[line width=1pt,color=colori](0.25,4.33)--(0.5,3.38)--(0.75,3.16)--(1.0,3.05)--(1.25,2.99)--(1.5,2.95)--(1.75,2.92)--(2.0,2.9)--(2.25,2.88)--(2.5,2.87)--(2.75,2.85)--(3.0,2.85)--(3.25,2.84)--(3.5,2.83)--(3.75,2.83)--(4.0,2.82)--(4.25,2.82)--(4.5,2.81)--(4.75,2.81)--(5.0,2.81)--(5.25,2.8)--(5.5,2.8)--(5.75,2.8)--(6.0,2.8)--(6.25,2.8)--(6.5,2.79)--(6.75,2.79)--(7.0,2.79)--(7.25,2.79)--(7.5,2.79)--(7.75,2.79)--(8.0,2.79)--(8.25,2.78)--(8.5,2.78); 
            \draw[line width=1pt,color=colori](6.4,5.0)--(7.2,5.0);\node at(7.2,5.0)[right]{$f(2n,n)$}; 
            \draw[line width=1pt,color=colorii](0.25,4.18)--(0.5,3.29)--(0.75,3.11)--(1.0,3.01)--(1.25,2.96)--(1.5,2.92)--(1.75,2.9)--(2.0,2.88)--(2.25,2.86)--(2.5,2.85)--(2.75,2.84)--(3.0,2.83)--(3.25,2.83)--(3.5,2.82)--(3.75,2.82)--(4.0,2.81)--(4.25,2.81)--(4.5,2.81)--(4.75,2.8)--(5.0,2.8)--(5.25,2.8)--(5.5,2.8)--(5.75,2.79)--(6.0,2.79)--(6.25,2.79)--(6.5,2.79)--(6.75,2.79)--(7.0,2.79)--(7.25,2.78)--(7.5,2.78)--(7.75,2.78)--(8.0,2.78)--(8.25,2.78); 
            \draw[line width=1pt,color=colorii](6.4,4.5)--(7.2,4.5);\node at(7.2,4.5)[right]{$f(3n,n)$}; 
            \draw[line width=1pt,densely dotted](0.25,0.15)--(0.5,0.62)--(0.75,1.1)--(1.0,1.39)--(1.25,1.59)--(1.5,1.75)--(1.75,1.86)--(2.0,1.96)--(2.25,2.03)--(2.5,2.09)--(2.75,2.15)--(3.0,2.19)--(3.25,2.23)--(3.5,2.26)--(3.75,2.29)--(4.0,2.32)--(4.25,2.34)--(4.5,2.36)--(4.75,2.38)--(5.0,2.4)--(5.25,2.41)--(5.5,2.43)--(5.75,2.44)--(6.0,2.45)--(6.25,2.46)--(6.5,2.47)--(6.75,2.48)--(7.0,2.49)--(7.25,2.5)--(7.5,2.51)--(7.75,2.52)--(8.0,2.52)--(8.25,2.53)--(8.5,2.53); 
            \draw[line width=1pt,densely dotted](6.4,1.0)--(7.2,1.0);\node at(7.2,1.0)[right]{$f^-(n,n)$}; 
            \draw[line width=1pt,densely dotted,color=colori](0.25,0.26)--(0.5,0.97)--(0.75,1.4)--(1.0,1.66)--(1.25,1.83)--(1.5,1.96)--(1.75,2.06)--(2.0,2.13)--(2.25,2.19)--(2.5,2.24)--(2.75,2.28)--(3.0,2.32)--(3.25,2.35)--(3.5,2.38)--(3.75,2.4)--(4.0,2.42)--(4.25,2.44)--(4.5,2.45)--(4.75,2.47)--(5.0,2.48)--(5.25,2.49)--(5.5,2.51)--(5.75,2.52)--(6.0,2.52)--(6.25,2.53)--(6.5,2.54)--(6.75,2.55)--(7.0,2.56)--(7.25,2.56)--(7.5,2.57)--(7.75,2.57)--(8.0,2.58)--(8.25,2.58)--(8.5,2.58); 
            \draw[line width=1pt,densely dotted,color=colori](6.4,1.5)--(7.2,1.5);\node at(7.2,1.5)[right]{$f^-(2n,n)$}; 
            \draw[line width=1pt,densely dotted,color=colorii](0.25,0.46)--(0.5,1.12)--(0.75,1.52)--(1.0,1.76)--(1.25,1.92)--(1.5,2.04)--(1.75,2.13)--(2.0,2.2)--(2.25,2.25)--(2.5,2.3)--(2.75,2.33)--(3.0,2.36)--(3.25,2.39)--(3.5,2.42)--(3.75,2.44)--(4.0,2.45)--(4.25,2.47)--(4.5,2.49)--(4.75,2.5)--(5.0,2.51)--(5.25,2.52)--(5.5,2.53)--(5.75,2.54)--(6.0,2.55)--(6.25,2.56)--(6.5,2.56)--(6.75,2.57)--(7.0,2.58)--(7.25,2.58)--(7.5,2.59)--(7.75,2.59)--(8.0,2.6)--(8.25,2.6); 
            \draw[line width=1pt,densely dotted,color=colorii](6.4,2.0)--(7.2,2.0);\node at(7.2,2.0)[right]{$f^-(3n,n)$}; 
        \end{tikzpicture}
        \caption{\small $f(m,n), f^-(m,n)$ for some fixed $\frac{m}{n}$}
         \label{fig:f_mn3}
    \end{figure}

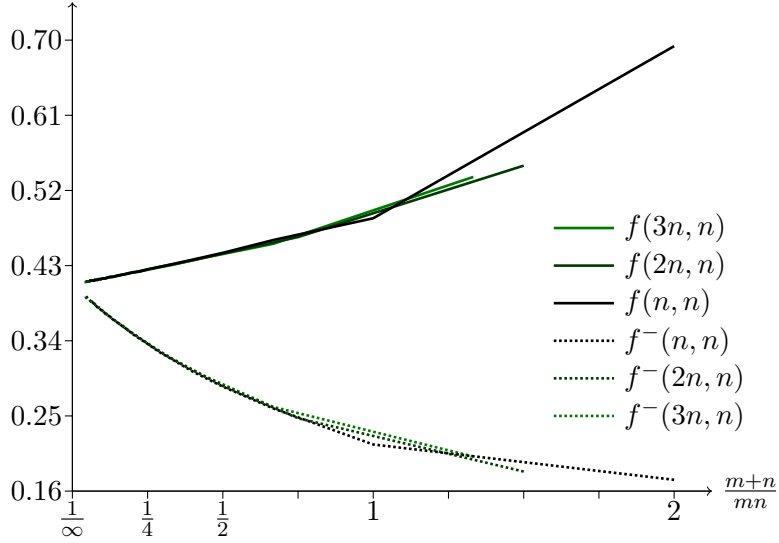
\begin{figure}
         \centering
         \begin{tikzpicture}[scale=1]
            \draw[->](0,0)--(8.5,0);\draw(0,-0.1)grid(8.5,0);\node at(0,0)[below]{$\frac{1}{\infty}$};\node at(1,0)[below]{$\frac{1}{4}$};\node at(2,0)[below]{$\frac{1}{2}$};\node at(4,0)[below]{1};\node at(8,0)[below]{2};\node at(8.5,0)[right]{$\frac{m+n}{mn}$}; 
            \draw[->](0,0)--(0,6.5);\draw(-0.1,0)grid(0,6.5);\node at(0,0)[left]{0.16};\node at(0,1)[left]{0.25};\node at(0,2)[left]{0.34};\node at(0,3)[left]{0.43};\node at(0,4)[left]{0.52};\node at(0,5)[left]{0.61};\node at(0,6)[left]{0.70};\node at(0, 6.5)[left]{}; 
            \draw[line width=1pt,color=colorii](5.33,4.18)--(2.67,3.29)--(1.78,3.11)--(1.33,3.01)--(1.07,2.96)--(0.89,2.92)--(0.76,2.9)--(0.67,2.88)--(0.59,2.86)--(0.53,2.85)--(0.48,2.84)--(0.44,2.83)--(0.41,2.83)--(0.38,2.82)--(0.36,2.82)--(0.33,2.81)--(0.31,2.81)--(0.3,2.81)--(0.28,2.8)--(0.27,2.8)--(0.25,2.8)--(0.24,2.8)--(0.23,2.79)--(0.22,2.79)--(0.21,2.79)--(0.21,2.79)--(0.2,2.79)--(0.19,2.79)--(0.18,2.78)--(0.18,2.78)--(0.17,2.78)--(0.17,2.78)--(0.16,2.78); 
            \draw[line width=1pt,color=colorii](6.4,3.5)--(7.2,3.5);\node at(7.2,3.5)[right]{$f(3n,n)$}; 
            \draw[line width=1pt,color=colori](6.0,4.33)--(3.0,3.38)--(2.0,3.16)--(1.5,3.05)--(1.2,2.99)--(1.0,2.95)--(0.86,2.92)--(0.75,2.9)--(0.67,2.88)--(0.6,2.87)--(0.55,2.85)--(0.5,2.85)--(0.46,2.84)--(0.43,2.83)--(0.4,2.83)--(0.38,2.82)--(0.35,2.82)--(0.33,2.81)--(0.32,2.81)--(0.3,2.81)--(0.29,2.8)--(0.27,2.8)--(0.26,2.8)--(0.25,2.8)--(0.24,2.8)--(0.23,2.79)--(0.22,2.79)--(0.21,2.79)--(0.21,2.79)--(0.2,2.79)--(0.19,2.79)--(0.19,2.79)--(0.18,2.78)--(0.18,2.78); 
            \draw[line width=1pt,color=colori](6.4,3.0)--(7.2,3.0);\node at(7.2,3.0)[right]{$f(2n,n)$}; 
            \draw[line width=1pt](8.0,5.92)--(4.0,3.63)--(2.67,3.34)--(2.0,3.17)--(1.6,3.08)--(1.33,3.02)--(1.14,2.98)--(1.0,2.95)--(0.89,2.92)--(0.8,2.91)--(0.73,2.89)--(0.67,2.88)--(0.62,2.87)--(0.57,2.86)--(0.53,2.85)--(0.5,2.85)--(0.47,2.84)--(0.44,2.83)--(0.42,2.83)--(0.4,2.83)--(0.38,2.82)--(0.36,2.82)--(0.35,2.82)--(0.33,2.81)--(0.32,2.81)--(0.31,2.81)--(0.3,2.81)--(0.29,2.8)--(0.28,2.8)--(0.27,2.8)--(0.26,2.8)--(0.25,2.8)--(0.24,2.8)--(0.24,2.79); 
            \draw[line width=1pt](6.4,2.5)--(7.2,2.5);\node at(7.2,2.5)[right]{$f(n,n)$}; 
            \draw[line width=1pt,densely dotted,color=colorii](5.33,0.46)--(2.67,1.12)--(1.78,1.52)--(1.33,1.76)--(1.07,1.92)--(0.89,2.04)--(0.76,2.13)--(0.67,2.2)--(0.59,2.25)--(0.53,2.3)--(0.48,2.33)--(0.44,2.36)--(0.41,2.39)--(0.38,2.42)--(0.36,2.44)--(0.33,2.45)--(0.31,2.47)--(0.3,2.49)--(0.28,2.5)--(0.27,2.51)--(0.25,2.52)--(0.24,2.53)--(0.23,2.54)--(0.22,2.55)--(0.21,2.56)--(0.21,2.56)--(0.2,2.57)--(0.19,2.58)--(0.18,2.58)--(0.18,2.59)--(0.17,2.59)--(0.17,2.6)--(0.16,2.6); 
            \draw[line width=1pt,densely dotted,color=colorii](6.4,1.0)--(7.2,1.0);\node at(7.2,1.0)[right]{$f^-(3n,n)$}; 
            \draw[line width=1pt,densely dotted,color=colori](6.0,0.26)--(3.0,0.97)--(2.0,1.4)--(1.5,1.66)--(1.2,1.83)--(1.0,1.96)--(0.86,2.06)--(0.75,2.13)--(0.67,2.19)--(0.6,2.24)--(0.55,2.28)--(0.5,2.32)--(0.46,2.35)--(0.43,2.38)--(0.4,2.4)--(0.38,2.42)--(0.35,2.44)--(0.33,2.45)--(0.32,2.47)--(0.3,2.48)--(0.29,2.49)--(0.27,2.51)--(0.26,2.52)--(0.25,2.52)--(0.24,2.53)--(0.23,2.54)--(0.22,2.55)--(0.21,2.56)--(0.21,2.56)--(0.2,2.57)--(0.19,2.57)--(0.19,2.58)--(0.18,2.58)--(0.18,2.59); 
            \draw[line width=1pt,densely dotted,color=colori](6.4,1.5)--(7.2,1.5);\node at(7.2,1.5)[right]{$f^-(2n,n)$}; 
            \draw[line width=1pt,densely dotted](8.0,0.15)--(4.0,0.62)--(2.67,1.1)--(2.0,1.39)--(1.6,1.59)--(1.33,1.75)--(1.14,1.86)--(1.0,1.96)--(0.89,2.03)--(0.8,2.09)--(0.73,2.15)--(0.67,2.19)--(0.62,2.23)--(0.57,2.26)--(0.53,2.29)--(0.5,2.32)--(0.47,2.34)--(0.44,2.36)--(0.42,2.38)--(0.4,2.4)--(0.38,2.41)--(0.36,2.43)--(0.35,2.44)--(0.33,2.45)--(0.32,2.46)--(0.31,2.47)--(0.3,2.48)--(0.29,2.49)--(0.28,2.5)--(0.27,2.51)--(0.26,2.52)--(0.25,2.52)--(0.24,2.53)--(0.24,2.54); 
            \draw[line width=1pt,densely dotted](6.4,2.0)--(7.2,2.0);\node at(7.2,2.0)[right]{$f^-(n,n)$}; 
        \end{tikzpicture}
         \caption{\small $f(m,n), f^-(m,n)$ for some fixed $\frac{m}{n}$, where the horizontal axis is $\frac{m+n}{mn}$ (density of boundary)}
         \label{fig:f_mn4}
    \end{figure}

    Figure~\ref{fig:f_mn} to~\ref{fig:f_mn4} present the values of $f(m,n)$ and $f^-(m,n)$ derived from our results for some specific cases.

    It can be observed that the combinatorial entropy of the model is primarily influenced by two effect:

    (1) Surface effect: The combinatorial entropy is larger for smaller grid sizes.
    This is because an $m \times n$ grid graph can be viewed as partitioning an infinite $\infty \times \infty$ grid graph into an $m \times n$ region.
    In contrast, the infinite grid relaxes adjacency restrictions along the partition boundaries, thereby increasing the entropy contribution from vertices near the boundary.
    As $m$ and $n$ tend to infinity, this increment is approximately proportional to the \textit{density} of boundary (i.e., the \textit{surface} for 2-dimensional cases), $\frac{m+n}{mn}$.

    Note that the ``density'' here is understood as the density of the boundary when an infinite number of finite grid graphs are pieced together into an infinite grid graph.
    Therefore, only one of the left and right sides of boundary is calculated, and the same applies to the upper and lower sides.
    In addition, the area is always defined as $mn$, including the case of $k^-$.

    (2) Parity effect: Having an odd length or width leads to a \textit{slightly} higher entropy compared to even dimensions, and this effect is compounded when both dimensions are odd.
    This occurs because, along edges of odd length, vertex selections can utilize space more ``efficiently'' than in even-length cases.
    In a sense, this effect almost only affect the area on the corner, because the existence of any a vertex on the boundary is equivalent to ``reset'' the parity of the boundary.
    Therefore, the parity effect is significant only when both the width and the height are small enough.

    In summary, as $m, n \to \infty$, the combinatorial entropy satisfies the \textit{area law}~\cite{area} as shown below:
    \begin{subequations}\label{eq:app}
        \begin{align}
            f(m,n) &\sim f + k\cdot\frac{m+n}{mn}+o(\frac{m+n}{mn}), \\
            f^-(m,n) &\sim f + k^-\cdot\frac{m+n+2}{(m+1)(n+1)}+o(\frac{m+n+2}{(m+1)(n+1)}),
        \end{align}
    \end{subequations}
    where $k$, $k^-$ are the \textit{surface effect constants}, which can be defined as\footnote{Approximation~\eqref{eq:app} can also be used to define these constants, but our definition relies on stronger convergence.}:
    \begin{subequations}\label{eq:k_def}
        \begin{align}
            k   &\coloneq \lim_{m,n\rightarrow\infty} \frac{f(m+\Delta m,n+\Delta n)-f(m,n)}{\frac{m+\Delta m+n+\Delta n}{(m+\Delta m)(n+\Delta n)}-\frac{m+n}{mn}}\\
                &= \lim_{m,n\rightarrow\infty} \frac{f(m,n)-f}{\frac{m+n}{mn}-0}; \\
            k^- &\coloneq \lim_{m,n\rightarrow\infty} \frac{f^-(m+\Delta m,n+\Delta n)-f^-(m,n)}{\frac{m+\Delta m+n+\Delta n}{(m+\Delta m)(n+\Delta n)}-\frac{m+n}{mn}}\\
                &= \lim_{m,n\rightarrow\infty} \frac{f^-(m,n)-f}{\frac{m+n}{mn}-0},
        \end{align}
    \end{subequations}
    where $\Delta m, \Delta n$ are any fixed positive integers.
    These two constants are the slope of extension lines at $x=0$ of the solid and dashed lines in Figure~\ref{fig:f_mn4}, respectively.
    Meanwhile, $f$ is the intercepts, and $y=f$ is the asymptote of the lines in Figure~\ref{fig:f_mn3}.
    Note that $m$ is fixed in Figure~\ref{fig:f_mn} and~\ref{fig:f_mn2} rather than tending towards infinity, so the lines in it have different convergence trends for different $m$.

    We cannot provide a rigorous proof of the convergence of $k$ or $k^-$ here.
    In fact, Forchhammer and Justesen~\cite{entropy1} proved that
    \[
        \frac{f(m,n)-f}{\frac{m+n}{mn}-0} \text{~~($\rightarrow k$, $m, n\rightarrow\infty$)}
    \]
    is bounded (greater than $0$ and with an upper bound), but this does not guarantee the convergence.
    And based on observation of the results, due to the existence of parity effect,
    \[
        \frac{f(m+\Delta m,n+\Delta n)-f(m,n)}{\frac{m+\Delta m+n+\Delta n}{(m+\Delta m)(n+\Delta n)}-\frac{m+n}{mn}} \text{~~($\rightarrow k$, $m, n\rightarrow\infty$)}
    \]
    is not monotonic with respect to $m$ and $n$ when $\Delta m=\Delta n=1$.
    However, from the observation of existing data, the corresponding data in all the models involved in this article seem convergent, even with the interference of parity effect.

    \subsection{Constant estimates}

    We will only give the estimation process of the constants of the version whose area is defined as $mn$, because both share the same $f$, and $k^-$ can be given by $k$ under our assumption of convergence:
    take $m=n$, we have
    \begin{subequations}
        \begin{align}
            k   &= \lim_{n\rightarrow\infty}\frac{n}{2}(f(n,n)-f);\\
            k^- &= \lim_{n\rightarrow\infty}\frac{n}{2}(f^-(n,n)-f)
                = \lim_{n\rightarrow\infty}\frac{n}{2}(\frac{n^2}{(n+1)^2}f(n,n)-f).
        \end{align}
    \end{subequations}
    By combining the two equations, we can get
    \begin{align}
        k^- &= \lim_{n\rightarrow\infty}\frac{n}{2}(\frac{n^2}{(n+1)^2}f(n,n)-(f(n,n)-\frac{k}{2n}))\\
            &= \lim_{n\rightarrow\infty}(-\frac{n(2n+1)}{2(n+1)^2}f(n,n)+k) \\
            &= k-f
    \end{align}
    To compute the two constants in this approximation, we first consider the combinatorial entropy of an $m \times \infty$ grid graph, where $\infty$ denotes either unidirectional or bidirectional infinite extension (the two cases are equivalent).
    The definition is:
    \begin{subequations}\label{eq:f_inf_def}
        \begin{align}
            f(m,\infty) &\coloneq \lim_{n \to \infty} f(m,n)\\
                        &=\frac{1}{m}\lim_{n \to \infty} \frac{\ln N(m,n)}{n}.
        \end{align}
    \end{subequations}
    So we can also denote $f$ as $f(\infty, \infty)$.
    And for each fixed $m$, we definite the surface effect constant that ignoring the upper and lower boundaries:
    \begin{subequations}\label{eq:k_inf_def}
        \begin{align}
            k(m,\infty) &\coloneq \lim_{n \to \infty} \frac{f(m+\Delta m,\infty)-f(m, \infty)}{\frac{1}{m+\Delta m}-\frac{1}{m}}\\
                        &= \lim_{n \to \infty} \frac{f(m,\infty)-f}{\frac{1}{m}-0}.
        \end{align}
    \end{subequations}
    That is, only the left and right boundaries are considered, and the density is $\frac{1}{n}$.

    We will use the results of $m\leq34, n\leq100$ for numerical analysis in the following text.
    The reason we calculating these additional results is, the estimation error comes from the residual term $o\left(\frac{m+n}{mn}\right)$ in~\eqref{eq:app}, and its absolute value is approximately proportional to $\frac{1}{mn}$.
    Therefore, it is best to use $m, n$ that can make $\frac{1}{mn}$ smaller for estimation.
    Here the complexity of the algorithm increases exponentially with respect to width $m$, while it only increases polynomial with respect to height $n$, so the ``narrower'' cases (with smaller $m$ but bigger $n$) can save more computational costs.

    We take data from two adjacent points and calculate the intercept and slope of their corresponding line:
    \begin{align}\label{eq:app1}
        \tilde{f}(n,\infty)(m,n)  &\coloneq \frac{\ln N(m,n+1)}{n+1}+ k(m,\infty)(n)\cdot\frac{1}{n+1};\\
        \tilde{k}(n,\infty)(m,n) &\coloneq \frac{\frac{\ln N(m,n+1)}{(n+1)}-\frac{\ln N(m,n)}{n}}{\frac{1}{n+1}-\frac{1}{n}},
    \end{align}
    to approach $f(m,\infty)$ and $k(m,\infty)$, respectively.
    Similarly, we can define
    \begin{align}\label{eq:app2}
        \tilde{f}(\infty,\infty)(m)  &\coloneq \frac{\ln N(m,\infty)}{m+1}+ k(\infty,\infty)(m)\cdot\frac{1}{m+1};\\
        \tilde{k}(\infty,\infty)(m) &\coloneq \frac{\frac{\ln N(m,\infty)}{(m+1)}-\frac{\ln N(m,\infty)}{m}}{\frac{1}{m+1}-\frac{1}{m}},
    \end{align}
    to approach $f$ and $k$, respectively.

    Although $f(m,n)$ ($f(m,\infty)$) can also be directily used to approximate $f(m,\infty)$ ($f(\infty, \infty)$), the approximation given by the above form converges faster\footnote{Pavlov~\cite{prob} observed that the convergence speed of $\tilde{f}(\infty,\infty)(m)$ is faster than exponential convergence, while $f(m,\infty)$ is slower.}.
    According to our assumption, this is because the bias caused by surface effect has been eliminated.

    Due to the existence of parity effect, the estimation of the constants given in this way always approaches the exact value from both sides according to the parity of the selected points, so that we can give the upper and lower bounds of the constants.
    When $m$ is fixed, the trend of odd and even terms of $f(m,\infty)$ and $k(m,\infty)$ always approaches exponential convergence.
    Therefore, we can retain its odd or even terms, and take the limit of the exponential curve passing through three adjacent points to obtain a more accurate estimate.

\begin{figure}
         \centering
         \begin{tikzpicture}[scale=1]
            \draw[->](0,0)--(8.5,0);\draw(0,-0.1)grid(8.5,0);\node at(0,0)[below]{0};\node at(1,0)[below]{4};\node at(2,0)[below]{8};\node at(3,0)[below]{12};\node at(4,0)[below]{16};\node at(5,0)[below]{20};\node at(6,0)[below]{24};\node at(7,0)[below]{28};\node at(8,0)[below]{32};\node at(8.5,0)[right]{$m$}; 
            \draw[->](0,0)--(0,6.5);\draw(-0.1,0)grid(0,6.5);\node at(0,0)[left]{0.29};\node at(0,1)[left]{0.32};\node at(0,2)[left]{0.35};\node at(0,3)[left]{0.38};\node at(0,4)[left]{0.41};\node at(0,5)[left]{0.44};\node at(0,6)[left]{0.47};\node at(0, 6.5)[left]{}; 
            \draw[line width=1pt](0.25,6.37)--(0.5,5.02)--(0.75,4.66)--(1.0,4.48)--(1.25,4.36)--(1.5,4.29)--(1.75,4.24)--(2.0,4.2)--(2.25,4.16)--(2.5,4.14)--(2.75,4.12)--(3.0,4.1)--(3.25,4.09)--(3.5,4.08)--(3.75,4.07)--(4.0,4.06)--(4.25,4.05)--(4.5,4.04)--(4.75,4.03)--(5.0,4.03)--(5.25,4.02)--(5.5,4.02)--(5.75,4.01)--(6.0,4.01)--(6.25,4.01)--(6.5,4.0)--(6.75,4.0)--(7.0,4.0)--(7.25,3.99)--(7.5,3.99)--(7.75,3.99)--(8.0,3.99)--(8.25,3.98)--(8.5,3.98); 
            \draw[line width=1pt](6.4,5.5)--(7.2,5.5);\node at(7.2,5.5)[right]{$f(m,\infty)$}; 
            \draw[line width=1pt,densely dashed](0.25,6.37)--(0.5,3.07)--(0.75,2.74)--(1.0,2.36)--(1.25,2.19)--(1.5,2.06)--(1.75,1.97)--(2.0,1.91)--(2.25,1.86)--(2.5,1.81)--(2.75,1.78)--(3.0,1.75)--(3.25,1.73)--(3.5,1.71)--(3.75,1.69)--(4.0,1.67)--(4.25,1.66)--(4.5,1.65)--(4.75,1.64)--(5.0,1.63)--(5.25,1.62)--(5.5,1.61)--(5.75,1.6)--(6.0,1.6)--(6.25,1.59)--(6.5,1.58)--(6.75,1.58)--(7.0,1.57)--(7.25,1.57)--(7.5,1.57)--(7.75,1.56)--(8.0,1.56)--(8.25,1.55)--(8.5,1.55); 
            \draw[line width=1pt,densely dashed](6.4,5.0)--(7.2,5.0);\node at(7.2,5.0)[right]{$f_\Delta (m,\infty)$}; 
            \draw[line width=1pt,densely dotted](0.25,6.37)--(0.5,1.89)--(0.75,1.83)--(1.0,1.21)--(1.25,1.09)--(1.5,0.89)--(1.75,0.81)--(2.0,0.72)--(2.25,0.66)--(2.5,0.61)--(2.75,0.57)--(3.0,0.53)--(3.25,0.5)--(3.5,0.48)--(3.75,0.46)--(4.0,0.44)--(4.25,0.42)--(4.5,0.41)--(4.75,0.39)--(5.0,0.38)--(5.25,0.37)--(5.5,0.36)--(5.75,0.35)--(6.0,0.34)--(6.25,0.34)--(6.5,0.33)--(6.75,0.32)--(7.0,0.32)--(7.25,0.31)--(7.5,0.31)--(7.75,0.3)--(8.0,0.3)--(8.25,0.29)--(8.5,0.29); 
            \draw[line width=1pt,densely dotted](6.4,4.5)--(7.2,4.5);\node at(7.2,4.5)[right]{$f_K (m,\infty)$}; 
        \end{tikzpicture}
         \caption{\small $f(m,\infty)$, $f_\varDelta(m,\infty)$ and $f_K(m,\infty)$}
         \label{fig:f_mxinf}
    \end{figure}
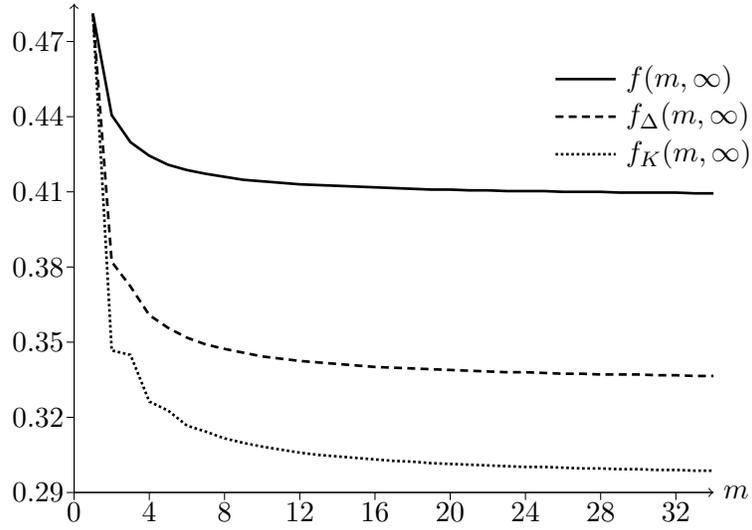

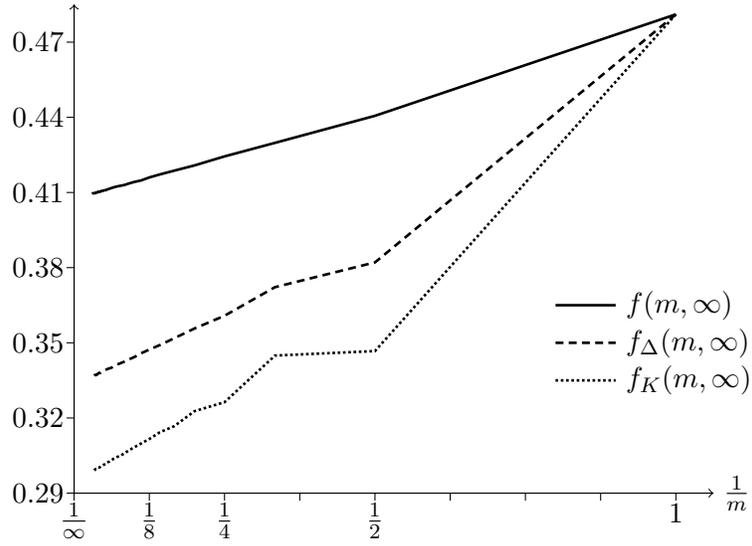
\begin{figure}
         \centering
         \begin{tikzpicture}[scale=1]
            \draw[->](0,0)--(8.5,0);\draw(0,-0.1)grid(8.5,0);\node at(0,0)[below]{$\frac{1}{\infty}$};\node at(1,0)[below]{$\frac{1}{8}$};\node at(2,0)[below]{$\frac{1}{4}$};\node at(4,0)[below]{$\frac{1}{2}$};\node at(8,0)[below]{1};\node at(8.5,0)[right]{$\frac{1}{m}$}; 
            \draw[->](0,0)--(0,6.5);\draw(-0.1,0)grid(0,6.5);\node at(0,0)[left]{0.29};\node at(0,1)[left]{0.32};\node at(0,2)[left]{0.35};\node at(0,3)[left]{0.38};\node at(0,4)[left]{0.41};\node at(0,5)[left]{0.44};\node at(0,6)[left]{0.47};\node at(0, 6.5)[left]{}; 
            \draw[line width=1pt](8.0,6.37)--(4.0,5.02)--(2.67,4.66)--(2.0,4.48)--(1.6,4.36)--(1.33,4.29)--(1.14,4.24)--(1.0,4.2)--(0.89,4.16)--(0.8,4.14)--(0.73,4.12)--(0.67,4.1)--(0.62,4.09)--(0.57,4.08)--(0.53,4.07)--(0.5,4.06)--(0.47,4.05)--(0.44,4.04)--(0.42,4.03)--(0.4,4.03)--(0.38,4.02)--(0.36,4.02)--(0.35,4.01)--(0.33,4.01)--(0.32,4.01)--(0.31,4.0)--(0.3,4.0)--(0.29,4.0)--(0.28,3.99)--(0.27,3.99)--(0.26,3.99)--(0.25,3.99)--(0.24,3.98)--(0.24,3.98); 
            \draw[line width=1pt](6.4,2.5)--(7.2,2.5);\node at(7.2,2.5)[right]{$f(m,\infty)$}; 
            \draw[line width=1pt,densely dashed](8.0,6.37)--(4.0,3.07)--(2.67,2.74)--(2.0,2.36)--(1.6,2.19)--(1.33,2.06)--(1.14,1.97)--(1.0,1.91)--(0.89,1.86)--(0.8,1.81)--(0.73,1.78)--(0.67,1.75)--(0.62,1.73)--(0.57,1.71)--(0.53,1.69)--(0.5,1.67)--(0.47,1.66)--(0.44,1.65)--(0.42,1.64)--(0.4,1.63)--(0.38,1.62)--(0.36,1.61)--(0.35,1.6)--(0.33,1.6)--(0.32,1.59)--(0.31,1.58)--(0.3,1.58)--(0.29,1.57)--(0.28,1.57)--(0.27,1.57)--(0.26,1.56)--(0.25,1.56)--(0.24,1.55)--(0.24,1.55); 
            \draw[line width=1pt,densely dashed](6.4,2.0)--(7.2,2.0);\node at(7.2,2.0)[right]{$f_\Delta (m,\infty)$}; 
            \draw[line width=1pt,densely dotted](8.0,6.37)--(4.0,1.89)--(2.67,1.83)--(2.0,1.21)--(1.6,1.09)--(1.33,0.89)--(1.14,0.81)--(1.0,0.72)--(0.89,0.66)--(0.8,0.61)--(0.73,0.57)--(0.67,0.53)--(0.62,0.5)--(0.57,0.48)--(0.53,0.46)--(0.5,0.44)--(0.47,0.42)--(0.44,0.41)--(0.42,0.39)--(0.4,0.38)--(0.38,0.37)--(0.36,0.36)--(0.35,0.35)--(0.33,0.34)--(0.32,0.34)--(0.31,0.33)--(0.3,0.32)--(0.29,0.32)--(0.28,0.31)--(0.27,0.31)--(0.26,0.3)--(0.25,0.3)--(0.24,0.29)--(0.24,0.29); 
            \draw[line width=1pt,densely dotted](6.4,1.5)--(7.2,1.5);\node at(7.2,1.5)[right]{$f_K (m,\infty)$}; 
        \end{tikzpicture}
         \caption{\small $f(m,\infty)$, $f_\varDelta(m,\infty)$ and $f_K(m,\infty)$, where the horizontal axis is $\frac{1}{m}$ (density of left and right boundary)}
         \label{fig:f_mxinf2}
    \end{figure}

    Figure~\ref{fig:f_mxinf} and Figure~\ref{fig:f_mxinf2} shows $f(m,\infty)$, where $m$ and $\frac{1}{m}$ as the horizontal axis, respectively.
    Figure~\ref{fig:k_mxinf} shows $k(m,\infty)$.
    More precise results are provided in Appendix B.

    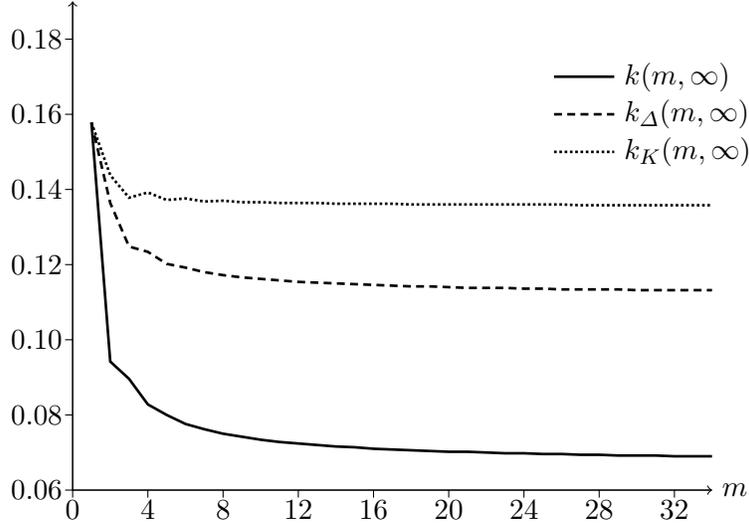
\begin{figure}
         \centering
         \begin{tikzpicture}[scale=1]
            \draw[->](0,0)--(8.5,0);\draw(0,-0.1)grid(8.5,0);\node at(0,0)[below]{0};\node at(1,0)[below]{4};\node at(2,0)[below]{8};\node at(3,0)[below]{12};\node at(4,0)[below]{16};\node at(5,0)[below]{20};\node at(6,0)[below]{24};\node at(7,0)[below]{28};\node at(8,0)[below]{32};\node at(8.5,0)[right]{$m$}; 
            \draw[->](0,0)--(0,6.5);\draw(-0.1,0)grid(0,6.5);\node at(0,0)[left]{0.06};\node at(0,1)[left]{0.08};\node at(0,2)[left]{0.10};\node at(0,3)[left]{0.12};\node at(0,4)[left]{0.14};\node at(0,5)[left]{0.16};\node at(0,6)[left]{0.18};\node at(0, 6.5)[left]{}; 
            \draw[line width=1pt](0.25,4.89)--(0.5,1.71)--(0.75,1.48)--(1.0,1.14)--(1.25,1.0)--(1.5,0.88)--(1.75,0.81)--(2.0,0.75)--(2.25,0.71)--(2.5,0.67)--(2.75,0.64)--(3.0,0.62)--(3.25,0.6)--(3.5,0.58)--(3.75,0.57)--(4.0,0.55)--(4.25,0.54)--(4.5,0.53)--(4.75,0.52)--(5.0,0.51)--(5.25,0.51)--(5.5,0.5)--(5.75,0.49)--(6.0,0.49)--(6.25,0.48)--(6.5,0.48)--(6.75,0.47)--(7.0,0.47)--(7.25,0.46)--(7.5,0.46)--(7.75,0.46)--(8.0,0.45)--(8.25,0.45)--(8.5,0.45); 
            \draw[line width=1pt](6.4,5.5)--(7.2,5.5);\node at(7.2,5.5)[right]{$k(m,\infty)$}; 
            \draw[line width=1pt,densely dashed](0.25,4.89)--(0.5,3.82)--(0.75,3.24)--(1.0,3.17)--(1.25,3.01)--(1.5,2.96)--(1.75,2.9)--(2.0,2.86)--(2.25,2.83)--(2.5,2.81)--(2.75,2.79)--(3.0,2.77)--(3.25,2.76)--(3.5,2.75)--(3.75,2.74)--(4.0,2.73)--(4.25,2.72)--(4.5,2.71)--(4.75,2.71)--(5.0,2.7)--(5.25,2.69)--(5.5,2.69)--(5.75,2.69)--(6.0,2.68)--(6.25,2.68)--(6.5,2.67)--(6.75,2.67)--(7.0,2.67)--(7.25,2.67)--(7.5,2.66)--(7.75,2.66)--(8.0,2.66)--(8.25,2.66)--(8.5,2.66); 
            \draw[line width=1pt,densely dashed](6.4,5.0)--(7.2,5.0);\node at(7.2,5.0)[right]{$k_\varDelta (m,\infty)$}; 
            \draw[line width=1pt,densely dotted](0.25,4.89)--(0.5,4.19)--(0.75,3.89)--(1.0,3.96)--(1.25,3.86)--(1.5,3.88)--(1.75,3.84)--(2.0,3.85)--(2.25,3.83)--(2.5,3.83)--(2.75,3.82)--(3.0,3.82)--(3.25,3.82)--(3.5,3.81)--(3.75,3.81)--(4.0,3.81)--(4.25,3.81)--(4.5,3.8)--(4.75,3.8)--(5.0,3.8)--(5.25,3.8)--(5.5,3.8)--(5.75,3.8)--(6.0,3.8)--(6.25,3.8)--(6.5,3.8)--(6.75,3.79)--(7.0,3.79)--(7.25,3.79)--(7.5,3.79)--(7.75,3.79)--(8.0,3.79)--(8.25,3.79)--(8.5,3.79); 
            \draw[line width=1pt,densely dotted](6.4,4.5)--(7.2,4.5);\node at(7.2,4.5)[right]{$k_K (m,\infty)$}; 
        \end{tikzpicture}
         \caption{\small $k(m,\infty)$, $k_\varDelta(m,\infty)$ and $k_K(m,\infty)$}
         \label{fig:k_mxinf}
    \end{figure}

    It can be seen that in the case of infinite in one direction, the parity effect in the other direction is not significant (but it still exists in reality).
    According to our previous description of parity effect, this can be understood as an infinitely long boundary increased the opportunity to reset parity in another direction.

    According to~\eqref{eq:app2}, we can obtain the upper and lower bounds of $f$ and $k$ by using the same method.
    However, the values of $f(m,\infty)$ and $k(m,\infty)$ we obtained already contain inherent errors, making direct estimation unreliable.
    To address this, we perform exponential interpolation using the values of $f$ and $k$ from three data points for each $m$: $(m, n)$, $(m, n-2)$, and $(m, n-4)$, and use the limit as the the approximate value of $f(m,\infty)$ and $k(m,\infty)$.
    Then, we increase $m$ and $n$ simultaneously and estimate the limit of the resulting sequence to obtain $f$ and $k$.
    \begin{align}
        0.407495101260688000450146812~3584 &< f<
        0.407495101260688000450146812~9046; \\
        1.5030480824753322643220663~2947 &< e^{f} <
        1.5030480824753322643220663~3030; \\
        0.0670552316597760242912072~6308 &< k<
        0.0670552316597760242912072~7191.
    \end{align}

    The OEIS sequence A085850 records a more accurate result of $e^{f}$ (which is called the \textit{entropy constant}), which is
    \begin{equation}
        e^{f}=1.5030480824753322643220663~294755536893857810\ldots.
    \end{equation}
    This shows that our estimate of $f$ is consistent with the results given by OEIS A085850, but not so accurate.
    However, there is no known results of the surface effect constant $k$ has been found before us.

    \section{Triangular grid graph and king graph}

    We denote the number of independent sets in a $n$-triangular grid graph as $N_\varDelta(n)$.
    Additionally, we denote the number of independent sets in a parallelogram-shaped triangular grid graph with width $m$ and oblique side length $n$ (we call this a $(m, n)$-triangular grid graph) as $N_\varDelta(m,n)$.
    An example of an independent set with $9$ vertices in a $(6, 4)$-triangular grid graph is shown below.
    \begin{center}
        \trigridii
    \end{center}

    An $n$-triangular grid graph can be embedded into the lower-right half of an $n \times n$ square grid graph. 
    Thus, it is equivalent to removing the upper-left half of the vertices from the $n \times n$ grid and adding edges connecting vertices that are adjacent along the anti-diagonal direction.

    By adding diagonal and anti-diagonal adjacency to the grid graph, we obtain the king graph (named after the movement pattern of chess kings)~\cite{king, king2}\footnote{Reference~\cite{king2} discussed the enumeration of tiling $1\times 1$ and $2\times 2$ squares in $m\times n$ rectangular area, which is obviously equivalent to $N_K(m-1,n-1)$.}.
    Consequently, this model is called the \textit{non-attacking kings} (NAK) model, since independent sets on the king graph correspond one-to-one with valid king placements on a chessboard where no kings attack each other.
    We denote the number of independent sets in a $m\times n$ king graph as $N_K(m,n)$.
    An example of an independent set with $5$ vertices in a $6 \times 4$ king graph is shown below.
    \begin{center}
        \kingi
    \end{center}

    \subsection{Algorithms}

    Following the approach in the previous section, for an independent set, we still select (\textbf{bold}) the edges above and to the left of each chosen vertex.
    The selected edges form non-overlapping L-shapes, which is the shortest W-shaped path, as the diagonal edges restrict the length of the paths. 

    \begin{center}
        \trigridi \raisebox{1.8cm}{$\Longrightarrow$} \trigridiii
    \end{center}

    To construct such paths, we simply remove template $\wtileiii{0}{0}{1}{1}$ from the original tensor, resulting in:
    \begin{equation}
        \varT_\varDelta\coloneq\wtileo{0}{0}{0}{0}+\wtilei{1}{0}{0}{0}+\wtileii{0}{1}{0}{1}+\wtileiv{0}{0}{1}{0}.
    \end{equation}

    Using this tensor, we assemble a triangular tensor network (that is, we remove the upper-right portion of the original tensor network above the main diagonal):
    \begin{equation}
        N_\varDelta (m,n)=\raisebox{-3cm}{\TNtri}
    \end{equation}
    and perform contractions to compute the independent set enumeration for triangular grid graphs.
    We have computed $N_\varDelta(n)$ for all $n \leq 40$, as shown in Appendix A, and also available in the corresponding OEIS sequence A027740.

    For computing $N_\varDelta(m,n)$, we can still use tensor $\varT_\varDelta$ but reverting to the rectangular tensor network shown in Figure 1.1.
    We have calculated all cases where $m + n \leq 77$, and the results are archived in OEIS sequences A219714 and A226444 (these two sequences essentially store the same dataset).

    And for $N_K(m,n)$, we modify the tensor configuration to a C-shaped path as illustrated below:
    \begin{equation}
        \kingii
    \end{equation}
    
    The corresponding tensor is given by:
    \begin{equation}
        \varT_K \coloneq\wtileo{0}{0}{0}{0}+\wtilei{1}{0}{0}{0}+\wtileii{0}{1}{0}{1}+\wtilev{0}{0}{1}{1}.
    \end{equation}
    
    We have also computed all cases where $m + n \leq 78$ (added the result of $m=n=39$ from A063443), with the results stored in OEIS sequence A245013.

    \subsection{Constant estimates}

    The variation trend of $N_\varDelta(m,n)$ and $N_K(m,n)$ are extremely close to $N(m,n)$, so we can use the method in the previous chapter to define (following Equation~\eqref{eq:f_def}\eqref{eq:k_def}\eqref{eq:f_inf_def}\eqref{eq:k_inf_def}) and estimate the combinatorial entropy and surface effect constant,
    and similarly perform numerical analysis using the results for all cases where $m\leq 35,n\leq 100$. 
    The algorithm is similar to that of $N(m,n)$, just replace the tensor $\varT$ in the Euqation~\eqref{eq:TN} with $\varT_\varDelta$ and $\varT_K$, obtained the upper and lower bounds:
    \begin{align}
        0.33324272197618188785374776~3953 &<f_\varDelta <
        0.33324272197618188785374776~4006;\\
        1.395485972479302735229500663~4939 &<e^{f_\varDelta}<
        1.395485972479302735229500663~5675;\\
        0.11182308857658~5958 &<k_\varDelta <
        0.11182308857658~6865;\\
        0.294640767816144918~2294 &<f_K <
        0.294640767816144918~4083;\\
        1.34264395112460129~7851 &<e^{f_K}<
        1.34264395112460129~8092;\\
        0.1354918026~4626 &<k_K <
        0.1354918026~7005.
    \end{align}
    and the precise value of $e^{f_\varDelta}$ (which can be obtained by solving the polynomial) is
    \begin{equation}
        e^{f_\varDelta}=1.395485972479302735229500663~566888\ldots
    \end{equation}
    which is consistent with our results.
    And the most accurate lower and upper bounds for $e^{f_K}$ are (given by ~\cite{lower, upper}):
    \begin{equation}
         1.34264395112460~1297\ldots <e^{f_K}
        <1.34264395112460~2238
    \end{equation}
    Our lower and upper bounds give three more exact figures compared to these.
    But due to the need for numerical analysis, which assumed that the data are sufficiently good, they are not as rigorous in mathematics.

    Likewise, Figure~\ref{fig:f_mxinf} and~\ref{fig:k_mxinf} show the corresponding constants for the one-sided infinite case.
    More precise results are given in Appendix B.
    Through comparison, it can be clearly seen that the change trend of the corresponding constants of the three kinds of graphs are roughly the same, but converges to different values.

    \section{Cylindrical grid graph}

    We denote the enumeration of independent sets on an $m\times n$ cylindrical grid graph (with periodic boundary conditions in the horizontal direction) as $N_C(m,n)$.

    The algorithms for cylindrical version is just to splicing the tensor network in the same way, and the tensor on the remaining boundaries remains unchanged.

    Following the approach in the previous chapter, we define the corresponding constants $f_{C}(m,n)$ (including the cases of $m=\infty$ or $n=\infty$).
    However, since the boundary conditions now differ in the two directions, the surface effect constants must also be direction-dependent.
    We define the surface effect constant corresponding to the periodic (horizontal) boundary as (assuming its convergence):
    \begin{subequations}\label{eq:k_cyc}
        \begin{align}
        \mathring{k}_{C} &\coloneq \lim_{m,n \to \infty} \frac{f_C(m,n+\Delta n) - f_C(m,n)}{\frac{m+n+\Delta n}{m(n+\Delta n)} - \frac{m+n}{mn}} \\
            &= \lim_{m,n \to \infty} \frac{f_C(m,n) - f_C(m,\infty)}{\frac{m+n}{mn} - 0},
    \end{align}
    \end{subequations}
    where \(\Delta n\) is an arbitrary fixed positive integer.
    And for the closed (vertical) boundary, it is defined analogously by swapping $m$ and $n$:
    \begin{subequations}
        \begin{align}
        \bar{k}_C &\coloneq \lim_{m,n \to \infty} \frac{f_C(m+\Delta m,n) - f_C(m,n)}{\frac{m+\Delta m+n}{(m+\Delta m)n+} - \frac{m+n}{mn}} \\
            &= \lim_{m,n \to \infty} \frac{f_C(m,n) - f_C(\infty,n)}{\frac{m+n}{mn} - 0},
    \end{align}
    \end{subequations}
    By these definitions, it can be easily shown that when this definition is applied to cases where $m$ and $n$ are interchangeable (as in all the previously discussed case), the resulting $\mathring{k}_C$ and $\bar{k}_C$ both coincide with the coefficient $k$ defined in Equation~\eqref{eq:k_def}.

    Note that the cylindrical grid graph of size $m\times n$ can also be assembled by splicing an $m\times n$ grid graph, so, similar to~\eqref{eq:k_def}, we can define:
    \begin{subequations}
        \begin{align}
            f_C(m,n)   &\coloneq \frac{\ln N_C (m,n)}{mn}; \\
            f^-_C(m,n) &\coloneq \frac{\ln N_C (m,n)}{(m+1)n}.
        \end{align}
    \end{subequations}
    Here the definition of area of the ``${}^-$'' version here is $(m+1)n$ rather than $(m+1)(n+1)$, because there is no splicing in the vertical direction.
    And we have:
    \begin{equation}
        f_C(m,n) < f_C <f^-_C(m,n);
    \end{equation}

    \begin{equation}
        \lim_{m,n\rightarrow\infty} f_C(m,n) =f_C \coloneq \lim_{m,n\rightarrow\infty} f^-_C(m,n)
    \end{equation}

    Substituted these into the Inequalities~\eqref{eq:ineq}, we get:
    \begin{equation}
        f^-(m,n)<f_C(m+1,n) <f(m+2,n),
    \end{equation}
    and,
    \begin{equation}
        f(\infty,n)=f_C(\infty,n)
    \end{equation}
    This can also be considered that when the width is infinite, the cylindrical and non-cylindrical grid graphs are equivalent.
    By definitions, these easily lead to:
    \begin{align}\label{eq:f_C}
        f_C       &=f, \\
        \bar{k}_C &=k.
    \end{align}

    Following the previous process, we calculate all cases with $m\leq31, n\leq91$ and perform numerical analysis.
    The range of our calculation here is slightly smaller than in the previous case, because the state tensor in the cylindrical case is larger.
    Some results of the cylindrical case are saved in the corresponding OEIS sequence A286513.
    Figure~\ref{fig:f_Cmxinf} to~\ref{fig:k_C2mxinf} show the values of $f_C (m,\infty)$ and $\bar{k}_C (m,\infty)$.

    \begin{figure}
         \centering
         \begin{tikzpicture}[scale=1]
            \draw[->](0,0)--(8.5,0);\draw(0,-0.1)grid(8.5,0);\node at(0,0)[below]{0};\node at(1,0)[below]{4};\node at(2,0)[below]{8};\node at(3,0)[below]{12};\node at(4,0)[below]{16};\node at(5,0)[below]{20};\node at(6,0)[below]{24};\node at(7,0)[below]{28};\node at(8,0)[below]{32};\node at(8.5,0)[right]{$m$}; 
            \draw[->](0,0)--(0,6.5);\draw(-0.1,0)grid(0,6.5);\node at(0,0)[left]{0.390};\node at(0,1)[left]{0.405};\node at(0,2)[left]{0.420};\node at(0,3)[left]{0.435};\node at(0,4)[left]{0.450};\node at(0,5)[left]{0.465};\node at(0,6)[left]{0.480};\node at(0, 6.5)[left]{}; 
            \draw[line width=1pt](0.25,6.08)--(0.5,3.38)--(0.75,2.66)--(1.0,2.28)--(1.25,2.06)--(1.5,1.91)--(1.75,1.8)--(2.0,1.73)--(2.25,1.66)--(2.5,1.61)--(2.75,1.57)--(3.0,1.54)--(3.25,1.51)--(3.5,1.49)--(3.75,1.46)--(4.0,1.45)--(4.25,1.43)--(4.5,1.41)--(4.75,1.4)--(5.0,1.39)--(5.25,1.38)--(5.5,1.37)--(5.75,1.36)--(6.0,1.35)--(6.25,1.35)--(6.5,1.34)--(6.75,1.33)--(7.0,1.33)--(7.25,1.32)--(7.5,1.32)--(7.75,1.31); 
            \draw[line width=1pt](6.4,5.5)--(7.2,5.5);\node at(7.2,5.5)[right]{$f(m,\infty)$}; 
            \draw[line width=1pt,densely dashed](0.5,3.38)--(0.75,0.55)--(1.0,1.34)--(1.25,1.11)--(1.5,1.19)--(1.75,1.16)--(2.0,1.17)--(2.25,1.17)--(2.5,1.17)--(2.75,1.17)--(3.0,1.17)--(3.25,1.17)--(3.5,1.17)--(3.75,1.17)--(4.0,1.17)--(4.25,1.17)--(4.5,1.17)--(4.75,1.17)--(5.0,1.17)--(5.25,1.17)--(5.5,1.17)--(5.75,1.17)--(6.0,1.17)--(6.25,1.17)--(6.5,1.17)--(6.75,1.17)--(7.0,1.17)--(7.25,1.17)--(7.5,1.17)--(7.75,1.17); 
            \draw[line width=1pt,densely dashed](6.4,5.0)--(7.2,5.0);\node at(7.2,5.0)[right]{$f_C(m,\infty)$}; 
         \end{tikzpicture}
         \caption{\small Comparison of $f(m,\infty)$ and $f_{C}(m,\infty)$}
         \label{fig:f_Cmxinf}
    \end{figure}
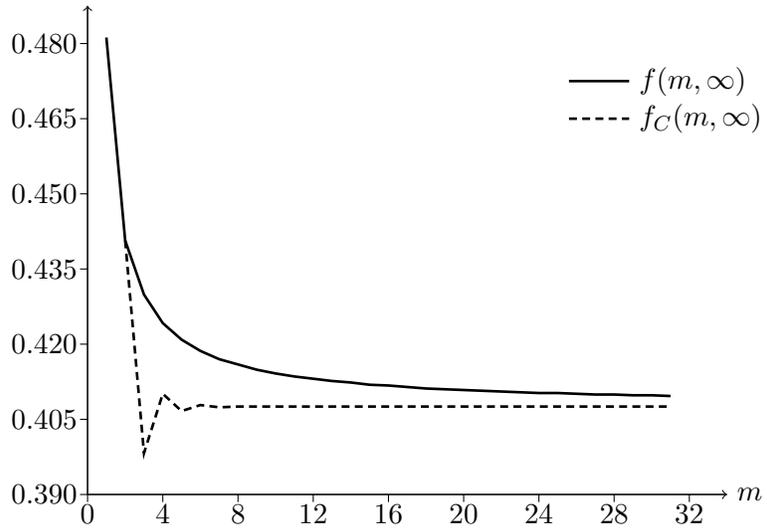

    \begin{figure}
         \centering
         \begin{tikzpicture}[scale=1]
            \draw[->](0,0)--(8.5,0);\draw(0,-0.1)grid(8.5,0);\node at(0,0)[below]{$\frac{1}{\infty}$};\node at(1,0)[below]{$\frac{1}{8}$};\node at(2,0)[below]{$\frac{1}{4}$};\node at(4,0)[below]{$\frac{1}{2}$};\node at(8,0)[below]{1};\node at(8.5,0)[right]{$\frac{1}{m}$}; 
            \draw[->](0,0)--(0,6.5);\draw(-0.1,0)grid(0,6.5);\node at(0,0)[left]{0.390};\node at(0,1)[left]{0.405};\node at(0,2)[left]{0.420};\node at(0,3)[left]{0.435};\node at(0,4)[left]{0.450};\node at(0,5)[left]{0.465};\node at(0,6)[left]{0.480};\node at(0, 6.5)[left]{}; 
            \draw[line width=1pt](8.0,6.08)--(4.0,3.38)--(2.67,2.66)--(2.0,2.28)--(1.6,2.06)--(1.33,1.91)--(1.14,1.8)--(1.0,1.73)--(0.89,1.66)--(0.8,1.61)--(0.73,1.57)--(0.67,1.54)--(0.62,1.51)--(0.57,1.49)--(0.53,1.46)--(0.5,1.45)--(0.47,1.43)--(0.44,1.41)--(0.42,1.4)--(0.4,1.39)--(0.38,1.38)--(0.36,1.37)--(0.35,1.36)--(0.33,1.35)--(0.32,1.35)--(0.31,1.34)--(0.3,1.33)--(0.29,1.33)--(0.28,1.32)--(0.27,1.32)--(0.26,1.31); 
            \draw[line width=1pt](6.4,3.5)--(7.2,3.5);\node at(7.2,3.5)[right]{$f(m,\infty)$}; 
            \draw[line width=1pt,densely dashed](4.0,3.38)--(2.67,0.55)--(2.0,1.34)--(1.6,1.11)--(1.33,1.19)--(1.14,1.16)--(1.0,1.17)--(0.89,1.17)--(0.8,1.17)--(0.73,1.17)--(0.67,1.17)--(0.62,1.17)--(0.57,1.17)--(0.53,1.17)--(0.5,1.17)--(0.47,1.17)--(0.44,1.17)--(0.42,1.17)--(0.4,1.17)--(0.38,1.17)--(0.36,1.17)--(0.35,1.17)--(0.33,1.17)--(0.32,1.17)--(0.31,1.17)--(0.3,1.17)--(0.29,1.17)--(0.28,1.17)--(0.27,1.17)--(0.26,1.17); 
            \draw[line width=1pt,densely dashed](6.4,3.0)--(7.2,3.0);\node at(7.2,3.0)[right]{$f_C(m,\infty)$}; 
         \end{tikzpicture}
         \caption{\small Comparison of $f(m,\infty)$ and $f_{C}(m,\infty)$, where the horizontal axis is $\frac{1}{m}$ (density of left and right boundary)}
         \label{fig:f_Cmxinf2}
    \end{figure}
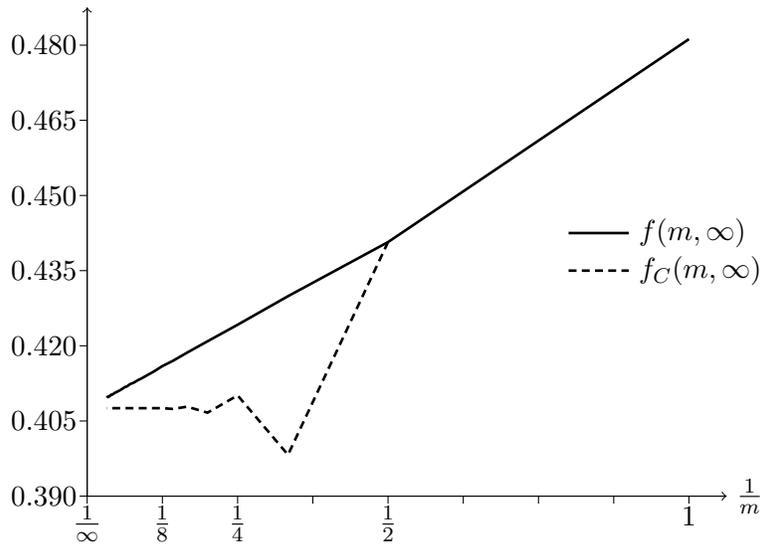

    \begin{figure}
         \centering
         \begin{tikzpicture}[scale=1]
            \draw[->](0,0)--(8.5,0);\draw(0,-0.1)grid(8.5,0);\node at(0,0)[below]{0};\node at(1,0)[below]{4};\node at(2,0)[below]{8};\node at(3,0)[below]{12};\node at(4,0)[below]{16};\node at(5,0)[below]{20};\node at(6,0)[below]{24};\node at(7,0)[below]{28};\node at(8,0)[below]{32};\node at(8.5,0)[right]{$m$}; 
            \draw[->](0,0)--(0,6.5);\draw(-0.1,0)grid(0,6.5);\node at(0,0)[left]{0.04};\node at(0,1)[left]{0.06};\node at(0,2)[left]{0.08};\node at(0,3)[left]{0.10};\node at(0,4)[left]{0.12};\node at(0,5)[left]{0.14};\node at(0,6)[left]{0.16};\node at(0, 6.5)[left]{}; 
            \draw[line width=1pt](0.25,5.89)--(0.5,2.71)--(0.75,2.48)--(1.0,2.14)--(1.25,2.0)--(1.5,1.88)--(1.75,1.81)--(2.0,1.75)--(2.25,1.71)--(2.5,1.67)--(2.75,1.64)--(3.0,1.62)--(3.25,1.6)--(3.5,1.58)--(3.75,1.57)--(4.0,1.55)--(4.25,1.54)--(4.5,1.53)--(4.75,1.52)--(5.0,1.51)--(5.25,1.51)--(5.5,1.5)--(5.75,1.49)--(6.0,1.49)--(6.25,1.48)--(6.5,1.48)--(6.75,1.47)--(7.0,1.47)--(7.25,1.46)--(7.5,1.46)--(7.75,1.46); 
            \draw[line width=1pt](6.4,5.5)--(7.2,5.5);\node at(7.2,5.5)[right]{$k(m,\infty)$}; 
            \draw[line width=1pt,densely dashed](0.5,2.71)--(0.75,0.95)--(1.0,1.47)--(1.25,1.32)--(1.5,1.36)--(1.75,1.35)--(2.0,1.35)--(2.25,1.35)--(2.5,1.35)--(2.75,1.35)--(3.0,1.35)--(3.25,1.35)--(3.5,1.35)--(3.75,1.35)--(4.0,1.35)--(4.25,1.35)--(4.5,1.35)--(4.75,1.35)--(5.0,1.35)--(5.25,1.35)--(5.5,1.35)--(5.75,1.35)--(6.0,1.35)--(6.25,1.35)--(6.5,1.35)--(6.75,1.35)--(7.0,1.35)--(7.25,1.35)--(7.5,1.35)--(7.75,1.35); 
            \draw[line width=1pt,densely dashed](6.4,5.0)--(7.2,5.0);\node at(7.2,5.0)[right]{$\bar{k}_C(m,\infty)$}; 
         \end{tikzpicture}
         \caption{\small Comparison of $k(m,\infty)$ and $\bar{k}_C(m,\infty)$}
         \label{fig:k_C2mxinf}
    \end{figure}
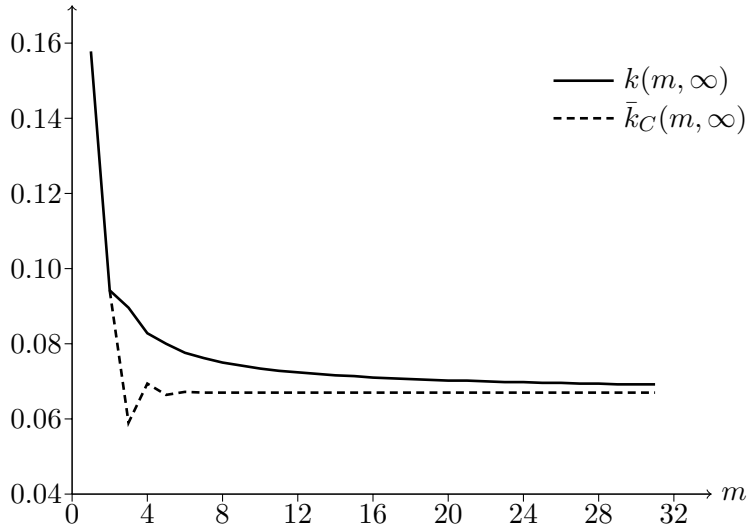

    Observations from the data reveal that the trend of $f_C(m,n)$ closely resembles that of $f(m,n)$, except that the horizontal boundary seems to exhibit no surface effect (can be observed from Figure~\ref{fig:f_Cmxinf} or~\ref{fig:f_Cmxinf2}, the asymptotic line is horizontal), i.e,
    \begin{equation}
        \mathring{k}_{C}=0,
    \end{equation}
    and the phase of the horizontal parity effect differs.

    Here we can use
    \begin{equation}
        \frac{f(m,\infty)-f(m-1,\infty)}{\frac{1}{m}-\frac{1}{m-1}} \text{~~($\rightarrow \mathring{k}_C, m\rightarrow\infty$)}
    \end{equation}
    to approach $\mathring{k}_C$.
    The following are the upper and lower bounds of the constants given by the results:
    \begin{align}
        0.407495101260~6521 &< f_C <
        0.407495101260~7636; \\
        1.503048082475~2784 &< e^{f_C} <
        1.503048082475~4459; \\
        -0.00000000000~1581 &< \mathring{k}_C <
         0.00000000000~3415; \\
        0.0670552316597~6690 &< \bar{k}_C <
        0.0670552316597~7922.
    \end{align}

    Since the obtained constants coincide with those in the planar version, we only want to know if evaluating these constants using periodic boundaries yields better precision.
    Clearly, these estimates are no more precise than those derived from the planar version.
    While this is partly due to the smaller data range, even with the same range of data, such estimations fail to enhance accuracy.
    This demonstrates that replacing boundaries with periodic boundaries only eliminates the first-order bias caused by surface effect, with negligible impact on higher-order bias (which determine the precision of such estimations).

    Using the similar method, it can be calculated and observed that the same conclusions hold true in the cylindrical versions of $(m,n)$-triangular grid graphs and king graphs (in this paper we omit the detailed discussion).
    However, strict proofs of these conclusions seem difficult to provide.

    \section{Summary and Further Expansions}

    In this paper, we present a method for expressing the independent set enumeration of grid graphs, $n$- and $(m,n)$-triangular grid graphs, king graphs, and cylindrical grid graphs as contractions of two-dimensional tensor networks.
    We also introduce a horizontal $l$-merging method applicable to these tensor networks, which reduces the spatial cost of the algorithm, thereby enabling computations for larger-scale systems.

    From these results, we observe that the behavior of the combinatorial entropy for finite instances of these models follows a similar trend, primarily influenced by surface effect and parity effect.
    Specifically, compared to the infinite case, the finite case exhibits a first-order bias in combinatorial entropy due to boundary density and a second-order deviation due to the parity of the boundary length.
    Based on observation of the results, we conjecture that the coefficients of the surface effect under our proposed definition are convergent.

    Through numerical analysis of the results, we estimate the combinatorial entropy for the corresponding infinite cases, as well as the surface effect constants in each direction (see Table~\ref{tab:f}).

    For cylindrical grid graphs, we prove that the surface effect constant for the closed boundary equals that of the original planar grid graph.
    As for the periodic boundary, we conjecture that the surface effect constant is zero, including the cylindrical version of $(m,n)$-triangular grid graphs and king graphs.
    We have reason to conjecture that many analogous models with periodic boundaries lack surface effect, though the precise conditions for this remain to be explored.

    By redesigning the tensor representation, our algorithm can potentially be extended to handle the three-dimensional version, but the computational cost will obviously increase further.
    Even if optimization is done through sparse matrices (where the effect is better than $l$-merging in this case), we can only calculate the cases of $l\times m\times n$ where $l, m\leq 6$.
    Based on these results\footnote{The results corresponding to $l=m=n\leq 6$ are saved in OEIS sequence A292702}, we can provide estimates of the corresponding constants:
    \begin{align}
        0.3622~2260 &< f_{\text{3d}} < 0.3622~8469; \\
        1.436~5186 &< e^{f_{\text{3d}}} < 1.436~6079; \\
        0.04~4592 &< k_{\text{3d}} < 0.04~5172.
    \end{align}
    The accuracy is significantly lower than that of the 2-dimensional case.
    In fact, it is even difficult to obtain such estimates for the corresponding 3-dimensional versions of other graphs because the convergence speed is slower.
    Therefore, the effective application of the method in 3-dimensional cases in this article still needs further exploration.

    In addition, other similar lattice-based statistical models, such as the Ising model and spin glass models, can also be computed using analogous methods, with similar constants defined and estimated.
    Reference~\cite{const} provides numerous examples of such models along with estimates of their corresponding entropy constants.
    However, it is worth noting that the efficiency of our algorithm relies critically on the reduction of tensor indices after $l$-merging, a feature that may not necessarily generalize to other models.

    Furthermore, by modifying specific terms within the tensor, we can compute \textit{weighted} entropy under predefined rules.
    Such weighted entropy carries clear physical significance: the system can be viewed as a two-dimensional crystal where the weights reflect potential energy contributions from particles with varying degrees of adjacency—the higher the potential energy, the lower the weight.
    However, the modular enumeration method becomes inapplicable if weights are non-integer.
    In addition, if configurations excluded by steric rules acquire non-zero weights after weighting, the efficiency gains from $l$-merging are no longer achievable.

    \bibliographystyle{unsrt}
    \bibliography{main2}

    \newpage
    \appendix
    
    \section{Results of $N_\varDelta(n)$}
    
    \tiny
    \begin{center}
        \begin{tabular}{ll}
            \hline
            $n$ & $N_\varDelta(n)$ \\
            \hline
            1 & 2 \\
            2 & 4 \\
            3 & 14 \\
            4 & 60 \\
            5 & 384 \\
            6 & 3318 \\
            7 & 40638 \\
            8 & 689636 \\
            9 & 16383974 \\
            10 & 542420394 \\
            11 & 25075022590 \\
            12 & 1617185558560 \\
            13 & 145563089994148 \\
            14 & 18283036276489970 \\
            15 & 3204638749437865046 \\
            16 & 783848125594781710150 \\
            17 & 267554112823378352976752 \\
            18 & 127443077148875066680607128 \\
            19 & 84712319011734340415085039268 \\
            20 & 78578236089316041630019193434422 \\
            21 & 101714646181932428248143602163208738 \\
            22 & 183734270908002233538130445637069602550 \\
            23 & 463150704646224649764463441631672026117820 \\
            24 & 1629220855732511492546360589163007528597792818 \\
            25 & 7997662422640270345091847114662902040439259297230 \\
            26 & 54786259684899846896275667859221962961228132566989930 \\
            27 & 523727899753149125104090756892359851817431315918615464448 \\
            28 & 6986589768922932982460178301569366275362668129029032113312408 \\
            29 & 130061947420267580631034025397303863700114859605321589571656334012 \\
            30 & 3378786393205609780342910249978713018052642672701712967887864644547490 \\
            31 & 122488889131163566432769056089046373073543507850377415047294050165791213246 \\
            32 & 6196666604535412306015914499142744511353655840972017434447326315502681041276430 \\
            33 & 437466740912168684635280144032063074869975883499801881532911772844870242164856276474 \\
            34 & 43098030681722787528219169313523445964284147519200826607903974311372797812802875833039502 \\
            35 & 5925094571292797096369643622887789155867678583221533368305832727550485439749459806470418956568 \\
            36 & 1136733404268645896674115362009195933978738542430435519788898362676587037502096685634287527216101776 \\
            37 & 304331869314329106476376132792082263561735057542101942691083735979066663899507450187178239886474017129262 \\
            38 & 113700328659089051056235801531311215244450101592511271633433453796597634278716541353621840952593885649619259026 \\
            39 & 59279083675659980886171478078363244554922103142394511114540011340057177851844108442007401890743518018775086303890120 \\
            40 & 43128734868728759732177228362983646007805350143316692609628166393492392639713173807759187796449840239360082182539803733482 \\
            \hline
        \end{tabular}
    \end{center}
    \normalsize

    \section{Results of $m\times\infty$ Cases}

    The values are the average of upper and lower bounds, and the errors are the distance between the values and the upper or lower bounds.

    Grid graph:

    \tiny
    ~~

    \begin{tabular}{rll}
        \hline
        $m$  & $f(m,\infty)$ & error\\
        \hline
        1 & 0.481211825059603447497758913424368423135184334385.. & 1.747e-82 \\
        2 & 0.440686793509771512616304662489896154514080164130.. & 1.254e-148 \\
        3 & 0.429871029469374617950502914215407621751493746077.. & 1.357e-73 \\
        4 & 0.424257111068728554960881148611310343385770258391.. & 3.886e-85 \\
        5 & 0.420906307591180525090652822971338741852067557428.. & 3.992e-62 \\
        6 & 0.418670960739873749692906536419297297610417504662.. & 1.446e-70 \\
        7 & 0.417074421385321816780333443033043436788244381734.. & 8.010e-52 \\
        8 & 0.415877005130446830774849291812225853508277913251.. & 8.448e-59 \\
        9 & 0.4149456825695727652037471297844666328527191144~2.. & 5.682e-47 \\
        10 & 0.414200624426308363715217595712867744428405846193.. & 2.674e-53 \\
        11 & 0.4135910314117670017779892681925885053004945~0484.. & 3.020e-44 \\
        12 & 0.413083037232346002699682398392625715616382528413.. & 3.174e-49 \\
        13 & 0.41265319600375189977128148105026539246011~4976.. & 1.440e-42 \\
        14 & 0.412284760664958275459409226752805875605855128~52.. & 1.451e-46 \\
        15 & 0.4119654500380065308691277760198135177550~7070.. & 1.775e-41 \\
        16 & 0.41168605323942402459214147892923319517674846~246.. & 9.756e-45 \\
        17 & 0.4114395266524395359156223124434991281951~2778.. & 9.592e-41 \\
        18 & 0.411220391908453336021773745187698496657393~1669.. & 1.952e-43 \\
        19 & 0.411024323979623580854538130414925550812~5065.. & 3.044e-40 \\
        20 & 0.41084786284367680170021684745069597390548~8288.. & 1.740e-42 \\
        21 & 0.410688207530201144471188690913678349738~1433.. & 6.731e-40 \\
        22 & 0.41054306633613236519210888233021034550034~8147.. & 8.851e-42 \\
        23 & 0.41041054611546087107180008510273441637~9154.. & 1.154e-39 \\
        24 & 0.4102890692465120014623411967975378700895~1127.. & 3.013e-41 \\
        25 & 0.41017731052707904142180761867366231500~1098.. & 1.639e-39 \\
        26 & 0.4100741486322178475382727582175709542277~9113.. & 7.633e-41 \\
        27 & 0.40997862835919822357204391589770528677~1656.. & 2.021e-39 \\
        28 & 0.409889930962822858460547184220356642248~3511.. & 1.549e-40 \\
        29 & 0.40980735062826648404639536242542927385~1737.. & 2.236e-39 \\
        30 & 0.409730275649347201259853726325633945829~8986.. & 2.651e-40 \\
        31 & 0.40965817324971303349179866085301097279~0761.. & 2.277e-39 \\
        32 & 0.409590577250056001209247039801263134634~3196.. & 3.972e-40 \\
        33 & 0.40952707797765091027715309334898084721~6025.. & 2.181e-39 \\
        34 & 0.409467313956563765870476437990537859488~9704.. & 5.363e-40 \\
        \hline
    \end{tabular}

    \begin{tabular}{rll}
        \hline
        $m$  & $k(m,\infty)$ & error \\
        \hline
        1 & 0.157704693902156707695138160235643026507568~1627.. & 1.713e-43 \\
        2 & 0.094113203229798857907688601760807870476330096950.. & 1.567e-75 \\
        3 & 0.089590595253620610958483031556189783985~0744.. & 9.287e-40 \\
        4 & 0.082706110770076437664639784430276786124850714074.. & 3.307e-50 \\
        5 & 0.0799636657371019442359679068953730957~7545.. & 3.502e-38 \\
        6 & 0.0776828631495366903604814876696446424868260~3315.. & 1.220e-44 \\
        7 & 0.076210158016163396359992393664850959~9250.. & 1.350e-37 \\
        8 & 0.07504915109469166499659477811134949817217~7525.. & 1.282e-42 \\
        9 & 0.074167201449674311189331903680694397~1554.. & 2.395e-37 \\
        10 & 0.0734535917564779054479654972801765690803~7262.. & 3.897e-41 \\
        11 & 0.072872869750062939058794970820749162~0830.. & 2.990e-37 \\
        12 & 0.072387688911533138399641388409320558971~3772.. & 6.629e-40 \\
        13 & 0.071977652512398651078419678585251967~2498.. & 3.119e-37 \\
        14 & 0.07162598863572782535655686417696289485~7281.. & 3.317e-39 \\
        15 & 0.071321297261801034384744194035040107~1100.. & 2.948e-37 \\
        16 & 0.07105465743353447688935973036592269292~8693.. & 9.278e-39 \\
        17 & 0.070819401596965957635253329296005247~0737.. & 3.729e-37 \\
        18 & 0.0706102791416443792397957778616521266~2653.. & 1.840e-38 \\
        19 & 0.070423172188881679028281970140296235~3504.. & 2.987e-37 \\
        20 & 0.0702547748162389155593212775996855676~7132.. & 2.924e-38 \\
        21 & 0.070102415767128742205132772972708911~7743.. & 2.127e-37 \\
        22 & 0.0699639073343422798000502260540626556~2622.. & 3.988e-38 \\
        23 & 0.069837443202453931831080970574661624~5808.. & 1.232e-37 \\
        24 & 0.0697215177093714253647846816922009571~5198.. & 7.604e-38 \\
        25 & 0.0696148662726917412578263442939397~2269.. & 3.659e-35 \\
        26 & 0.0695164187852539187271538693626759707~6526.. & 8.090e-38 \\
        27 & 0.0694252637075554206289884555248855134~2813.. & 3.224e-38 \\
        28 & 0.0693406197053978042029012221213195930~5661.. & 7.901e-38 \\
        29 & 0.0692618132212655882751137793392571195~3206.. & 7.420e-38 \\
        30 & 0.0691882605024607765257384775377960059~2528.. & 7.091e-38 \\
        31 & 0.0691194531204780203570578461962257885~6314.. & 9.503e-38 \\
        32 & 0.0690549461998135617462906663562608670~7103.. & 5.823e-38 \\
        33 & 0.0689943487895172154034152171991675924~6264.. & 9.867e-38 \\
        34 & 0.0689373159327566231401134135441448806~1161.. & 4.318e-38 \\
        \hline
    \end{tabular}

    ~~
    \normalsize

    \newpage
    Triangular grid graph:

    \tiny
    ~~

    \begin{tabular}{rll}
        \hline
        $m$  & $f_\varDelta(m,\infty)$ & error\\
        \hline
        1 & 0.481211825059603447497758913424368423135184334385.. & 1.747e-82 \\
        2 & 0.38224508584003564132935849918485739375941641~915.. & 3.349e-45 \\
        3 & 0.372185143199076763266438568096490403101529933154.. & 2.819e-49 \\
        4 & 0.3607372594505253932356497619912931246474689208~0.. & 2.230e-47 \\
        5 & 0.35574260139400142667383872264604317138623503951~.. & 6.381e-48 \\
        6 & 0.351838510551468457327336830230539179770776192~18.. & 1.091e-46 \\
        7 & 0.3492304710973854521753183718551848532972884234~0.. & 8.923e-47 \\
        8 & 0.347216425333335764595019305431779792865100194~23.. & 2.694e-46 \\
        9 & 0.345668874361900112165909451684223000385211678~27.. & 9.082e-46 \\
        10 & 0.344424579960030349578661266569394297742112559~87.. & 4.340e-46 \\
        11 & 0.343408607760006355731375982848042383722602510~15.. & 8.011e-46 \\
        12 & 0.342561262063819766309302250075389279069756298~73.. & 8.093e-46 \\
        13 & 0.3418445151749295896587401049344373088696597~4372.. & 1.001e-44 \\
        14 & 0.34123007959877623106055112000280299170447600~018.. & 1.569e-45 \\
        15 & 0.34069759655121014213426871695280823289617100~637.. & 7.458e-45 \\
        16 & 0.34023166432283304403008842353019848544135116~989.. & 6.473e-45 \\
        17 & 0.3398205509541514536439138349684883168071509~4929.. & 1.284e-44 \\
        18 & 0.3394551157034229582635409277398703885069682~4781.. & 3.265e-44 \\
        19 & 0.3391281477194529652638295786705659963786269~8875.. & 3.405e-44 \\
        20 & 0.338833876394998084558089899709528124848592~1676.. & 4.486e-43 \\
        21 & 0.338567630959519888013826215503303318356141~1154.. & 2.462e-43 \\
        22 & 0.338325589637527454336559674309339012930534~2153.. & 6.331e-43 \\
        23 & 0.338104595392985880697959795767070906878520~3133.. & 2.382e-43 \\
        24 & 0.337902017333388166149907579935765553552723~4030.. & 7.134e-43 \\
        25 & 0.337715645519298739666219293958362935060086~3543.. & 7.564e-43 \\
        26 & 0.33754360999833945353622915781621265905764~5568.. & 1.430e-42 \\
        27 & 0.33738431784939546073229378761077942603728~5567.. & 2.436e-42 \\
        28 & 0.33723640371105762413099245829477556477482~3970.. & 3.248e-42 \\
        29 & 0.33709869054778915895037690262621621525631~5198.. & 7.484e-42 \\
        30 & 0.33697015826206781234582488981430152444752~1126.. & 7.663e-42 \\
        31 & 0.3368499183818783360750809032067347006364~7320.. & 2.019e-41 \\
        32 & 0.3367371934942001825786336298848355354715~8657.. & 1.781e-41 \\
        33 & 0.3366313004178966475414730549578153917682~9039.. & 4.406e-41 \\
        34 & 0.3365316363460814898609100974438130774421~4908.. & 3.951e-41 \\
        \hline
    \end{tabular}

    \begin{tabular}{rll}
        \hline
        $m$  & $k_\varDelta(m,\infty)$ & error \\
        \hline
        1 & 0.157704693902156707695138160235643026507568~1627.. & 1.713e-43 \\
        2 & 0.136318375996289851824900656377108616529558~6145.. & 2.942e-43 \\
        3 & 0.124779734547815277575328401242676768879648440946.. & 1.657e-49 \\
        4 & 0.12340128633760295349358663830804577211694463~471.. & 1.964e-45 \\
        5 & 0.120269670957925531231970554977389421840926019~90.. & 5.806e-46 \\
        6 & 0.11919010179930900780970428779547415049113528~862.. & 9.437e-45 \\
        7 & 0.11801051816905876235207899422347042425553638~001.. & 7.825e-45 \\
        8 & 0.1172855408621138041544067145812395978210482~2096.. & 2.109e-44 \\
        9 & 0.1166602842871648936642580154725193435103388~2812.. & 8.216e-44 \\
        10 & 0.1161834562799485167860292613123980064093115~9213.. & 3.813e-44 \\
        11 & 0.1157844771411471693625783131620706642506280~4907.. & 7.089e-44 \\
        12 & 0.1154553251984579630498524072687294586420934~7530.. & 7.098e-44 \\
        13 & 0.115175563804177121239701172335088404939898~0421.. & 9.557e-43 \\
        14 & 0.114936234376557618379573906366753118180169~5803.. & 1.379e-43 \\
        15 & 0.114728642009003165045569858243177810498211~4506.. & 6.466e-43 \\
        16 & 0.114547063171153367928106540098894313891609~7504.. & 5.706e-43 \\
        17 & 0.11438682262232187740068207289391332953340~4823.. & 1.124e-42 \\
        18 & 0.11424439544688626417346344141735610648952~4933.. & 2.899e-42 \\
        19 & 0.11411695732242368356953819768221990334652~7526.. & 3.003e-42 \\
        20 & 0.1140022642248757471196083330967720102360~6188.. & 4.195e-41 \\
        21 & 0.1138984938306209342425785042848056594528~5284.. & 2.437e-41 \\
        22 & 0.1138041572743470220380823763281233612010~2767.. & 5.128e-41 \\
        23 & 0.1137180238356884460676627683959065914150~5011.. & 2.094e-41 \\
        24 & 0.1136390682061608161437306020812945788215~8048.. & 6.209e-41 \\
        25 & 0.1135664290186695165669693571728736934614~6456.. & 6.673e-41 \\
        26 & 0.113499377464054630448090685602368756819~8238.. & 1.251e-40 \\
        27 & 0.113437292690131898901798507927305962124~4462.. & 2.156e-40 \\
        28 & 0.113379642543334679915539864493437664232~3510.. & 2.850e-40 \\
        29 & 0.113325968268576876549615889259872726360~5232.. & 6.645e-40 \\
        30 & 0.113275872278859448989700868142222172987~4146.. & 6.741e-40 \\
        31 & 0.11322900828845782068981439249934444526~8575.. & 1.797e-39 \\
        32 & 0.11318507329746395365313807764011883336~1099.. & 1.569e-39 \\
        33 & 0.11314380103319416877820268162892073932~2538.. & 3.921e-39 \\
        34 & 0.11310495654917658496510327263170924234~5723.. & 3.478e-39 \\
        \hline
    \end{tabular}

    ~~
    \normalsize

    \newpage
    King graph:

    \tiny
    ~~

    \begin{tabular}{rll}
    \hline
    $m$  & $f_K(m,\infty)$ & error\\
    \hline
    2 & 0.346573590279972654708616060729088284037750067180.. & 3.734e-60 \\
    3 & 0.344822371107929602747312223453205280967874552348.. & 7.738e-64 \\
    4 & 0.326427375567256173180814992954693672829363073224.. & 4.998e-59 \\
    5 & 0.322590235671261671112407057675421432742672798331.. & 4.962e-49 \\
    6 & 0.31684368489562135692756142586633249968440~9214.. & 5.731e-42 \\
    7 & 0.314165959899889025046980745330400949105~9228.. & 6.189e-40 \\
    8 & 0.3114990195482019434221915248023320790~7571.. & 3.803e-38 \\
    9 & 0.309731861671801325410197272997702785~7578.. & 4.100e-37 \\
    10 & 0.30817267872520856461928504486417883~7723.. & 2.986e-36 \\
    11 & 0.3069664500987039731202584595685515~1027.. & 1.192e-35 \\
    12 & 0.3059277764743654622161015888861767~1957.. & 3.673e-35 \\
    13 & 0.3050651423949307315923604078365386~5379.. & 8.607e-35 \\
    14 & 0.304317812429161035087236388633332~5819.. & 1.714e-34 \\
    15 & 0.303674016432559217929075395066162~1009.. & 2.944e-34 \\
    16 & 0.303108777914932795994920029492631~7978.. & 4.557e-34 \\
    17 & 0.302610986437883918305506063734775~5045.. & 5.630e-34 \\
    18 & 0.302168034294963190824420837696186~2155.. & 7.545e-34 \\
    19 & 0.301771943158146334285728371529948~1827.. & 9.544e-34 \\
    20 & 0.30141534405071055749212804672475~7006.. & 1.153e-33 \\
    21 & 0.30109276537755882768198405587117~9823.. & 1.342e-33 \\
    22 & 0.30079948262808392374428960790034~7567.. & 1.514e-33 \\
    23 & 0.30053171751372096993966298091248~0996.. & 1.666e-33 \\
    24 & 0.30028625870980276076740654405578~7246.. & 1.797e-33 \\
    25 & 0.30006044036918319570177605713898~2746.. & 1.907e-33 \\
    26 & 0.29985199077023651447586907929830~0292.. & 1.996e-33 \\
    27 & 0.29965898284402334475373759006792~6008.. & 2.067e-33 \\
    28 & 0.29947976071079056985488672198281~7979.. & 2.122e-33 \\
    29 & 0.29931289897207662385429068679590~6642.. & 2.163e-33 \\
    30 & 0.29915716122356639049348736571627~6558.. & 2.193e-33 \\
    31 & 0.29901147113568235303332948692525~0898.. & 2.212e-33 \\
    32 & 0.29887488664572720326331468708603~3725.. & 2.224e-33 \\
    33 & 0.29874658002024660807452080095424~1243.. & 2.228e-33 \\
    34 & 0.29862582083485532272779794290053~5526.. & 2.227e-33 \\
    \hline
    \end{tabular}

    \begin{tabular}{rll}
    \hline

    $m$  & $k_K(m,\infty)$ & error \\
    \hline
    2 & 0.1438410362258904637196095029968~7488.. & 3.876e-32 \\
    3 & 0.137838640156139641716674857452752~9000.. & 3.222e-34 \\
    4 & 0.139286411793217535746569204330~7032.. & 1.076e-31 \\
    5 & 0.1371086196888887613016872007235~1793.. & 5.367e-32 \\
    6 & 0.137623353042322015274233198752~2794.. & 3.289e-31 \\
    7 & 0.136873799811907816736099425478~2138.. & 3.593e-31 \\
    8 & 0.136940042230733570194158522702~5581.. & 7.169e-31 \\
    9 & 0.136650278122244833795354613951~4036.. & 8.663e-31 \\
    10 & 0.13660271944217935392440719189~7246.. & 1.192e-30 \\
    11 & 0.13646556672674830660656631147~8811.. & 1.396e-30 \\
    12 & 0.13640344859407072664562085628~4368.. & 1.660e-30 \\
    13 & 0.13632332228673005744121204039~3978.. & 1.860e-30 \\
    14 & 0.13626916918760208909325781862~3680.. & 2.066e-30 \\
    15 & 0.13621460058584173240431074022~5476.. & 2.235e-30 \\
    16 & 0.13617086073671872778512114531~9506.. & 2.391e-30 \\
    17 & 0.13613016622712294028479883566~4209.. & 2.524e-30 \\
    18 & 0.13609509323267944414403415581~6881.. & 2.644e-30 \\
    19 & 0.13606313663954589725046910812~9837.. & 2.747e-30 \\
    20 & 0.13603467664277665957257760660~2458.. & 2.841e-30 \\
    21 & 0.13600876982841219612416300067~0067.. & 2.925e-30 \\
    22 & 0.13598530036441792174132084042~1071.. & 3.001e-30 \\
    23 & 0.13596382879336400576524589799~9215.. & 3.072e-30 \\
    24 & 0.13594416894006850733159554103~4131.. & 3.137e-30 \\
    25 & 0.13592607016791124692466706213~5950.. & 3.199e-30 \\
    26 & 0.13590936972142297194875316394~6521.. & 3.257e-30 \\
    27 & 0.13589390315392282761368263338~5445.. & 3.312e-30 \\
    28 & 0.13587954300715461157614243988~4716.. & 3.364e-30 \\
    29 & 0.13586617234567900406967116660~6948.. & 3.414e-30 \\
    30 & 0.13585369351561877542626767890~8996.. & 3.461e-30 \\
    31 & 0.13584201953437621247611629885~8882.. & 3.506e-30 \\
    32 & 0.13583107530068056288331564637~8560.. & 3.549e-30 \\
    33 & 0.13582079428928758574965475463~2768.. & 3.590e-30 \\
    34 & 0.13581111807699019744299045146~7159.. & 3.629e-30 \\
    \hline
    \end{tabular}

    ~~
    \normalsize

    \newpage
    Cylindrical grid graph:

    \tiny
    ~~

    \begin{tabular}{rll}
        \hline
        $m$  & $f_C(m,\infty)$ & error\\
        \hline
        1 & 0 (the only independent set is $\varnothing$) & 0 \\
        2 & 0.440686793509771512616304662489896154514080164130.. & 1.251e-129 \\
        3 & 0.398254405762369768037310276173030174512054025051.. & 1.741e-174 \\
        4 & 0.410056037580579286727537638792521264389586605666.. & 2.139e-100 \\
        5 & 0.406614574980743081640389836200888915263352129819.. & 2.033e-147 \\
        6 & 0.407805573398536382991903920818913345252317984757.. & 1.342e-85 \\
        7 & 0.407377738577143243451221510847062650686604173598.. & 8.335e-98 \\
        8 & 0.407540536130730140731839318252510680906485802109.. & 5.210e-85 \\
        9 & 0.407476963681491300463330442585951099346721397808.. & 1.295e-86 \\
        10 & 0.407502471475487026435527618723306480058269230819.. & 8.362e-88 \\
        11 & 0.407492054221542551040916475126610265319801278556.. & 1.230e-85 \\
        12 & 0.407496377100376185384519327125170551732763696616.. & 7.995e-88 \\
        13 & 0.407494560965162970545768569350680042886730654216.. & 2.953e-88 \\
        14 & 0.407495332195537740906668410468130049793076448197.. & 9.131e-82 \\
        15 & 0.407495001749216172449207026695309846175670097667.. & 4.803e-88 \\
        16 & 0.407495144440752752741969695175287476323400224471.. & 1.609e-77 \\
        17 & 0.407495082408498643751305339031118807868850953304.. & 4.176e-85 \\
        18 & 0.407495109536104074790465042795975811149297052744.. & 1.283e-74 \\
        19 & 0.407495097610523111547294519582805426452560404542.. & 5.537e-81 \\
        20 & 0.407495102877685270603811366618664293414579011448.. & 1.450e-72 \\
        21 & 0.407495100541578208244421788140535724922803069607.. & 6.013e-78 \\
        22 & 0.407495101581617887126121171542424854050914030074.. & 4.497e-71 \\
        23 & 0.407495101117001368589640269493627537165536599329.. & 1.174e-75 \\
        24 & 0.407495101325207882306455345585423473932062789344.. & 5.734e-70 \\
        25 & 0.407495101231638639587932465274376676688491194070.. & 6.773e-74 \\
        26 & 0.407495101273799510554827159084124619267373898347.. & 3.898e-69 \\
        27 & 0.407495101254756563506359230289825949211319186169.. & 1.609e-72 \\
        28 & 0.407495101263376965740585753999796650834757405018.. & 1.678e-68 \\
        29 & 0.407495101259466585205452618291974879303603309567.. & 1.972e-71 \\
        30 & 0.407495101261243822574379032361490716480249123388.. & 5.152e-68 \\
        31 & 0.407495101260434633484484050423978390730211167124.. & 1.463e-70 \\
        \hline
    \end{tabular}

    \begin{tabular}{rll}
        \hline
        $m$  & $k_C(m,\infty)$ & error \\
        \hline
        1 & -- & -- \\
        2 & 0.094113203229798857907688601760807870476330096950.. & 1.353e-68 \\
        3 & 0.058928544290683157847990049517066671735428179528.. & 1.187e-91 \\
        4 & 0.069384135414422895736721965515435849723603530817.. & 1.050e-87 \\
        5 & 0.066349971797211572898678998542916885396729949589.. & 1.119e-86 \\
        6 & 0.067271283366836719331099523813977714627897082687.. & 1.566e-75 \\
        7 & 0.066990956574619982750578768491331448486517499317.. & 5.186e-87 \\
        8 & 0.067073487883750124065906408790421520135592315860.. & 4.958e-69 \\
        9 & 0.067050807539628474176556638796077477603984151538.. & 1.366e-83 \\
        10 & 0.067055875265203839591691111087583241250917952574.. & 2.574e-66 \\
        11 & 0.067055449697414623250926893674904138368831987764.. & 2.645e-76 \\
        12 & 0.067054939670252881058575824970385362679783027662.. & 5.759e-65 \\
        13 & 0.067055439843000970264267785735087829089945918016.. & 3.688e-72 \\
        14 & 0.067055107004958819081878953863174320845937647593.. & 3.152e-64 \\
        15 & 0.067055300559273677912863560564681415520385524124.. & 1.353e-69 \\
        16 & 0.067055195268804030622893454747952449219041017221.. & 8.415e-64 \\
        17 & 0.067055250342055383438047886787789777632224661100.. & 6.534e-68 \\
        18 & 0.067055222251013947991134430036066893905291095666.. & 1.509e-63 \\
        19 & 0.067055236334546333312907660484709804517693129973.. & 9.321e-67 \\
        20 & 0.067055229360003603776196466269453314968536147662.. & 2.135e-63 \\
        21 & 0.067055232782770017809606609542990260809317479646.. & 6.122e-66 \\
        22 & 0.067055231114536093633528215212248391815161374633.. & 2.609e-63 \\
        23 & 0.067055231923324108412205910254879476862377664872.. & 2.403e-65 \\
        24 & 0.067055231532836387141123141319643440912450410544.. & 2.898e-63 \\
        25 & 0.067055231720744515806079458731005688214044471962.. & 6.613e-65 \\
        26 & 0.067055231630560092651377513594770571366742143632.. & 3.009e-63 \\
        27 & 0.067055231673750133496731804564530192471753485104.. & 1.417e-64 \\
        28 & 0.067055231653102360943547136134394664015312947581.. & 8.365e-65 \\
        29 & 0.067055231662959172577426601984300852713008382333.. & 2.528e-64 \\
        30 & 0.067055231658259336594700238327482607097573034216.. & 7.931e-65 \\
        31 & 0.067055231660498060685913047479039789986546042007.. & 3.935e-64 \\
        \hline
    \end{tabular}

    ~~
    \normalsize

\end{document}